\documentclass{amsart}

\usepackage{latexsym,amssymb,lastpage}
\usepackage{graphicx,amsfonts,amsbsy}
\usepackage{times,mathptmx,bm,amsmath}
\usepackage{color}
\usepackage{tikz}
\usepackage{subcaption}
\usepackage{hyperref}
\usepackage{epstopdf}
\usepackage{amsthm}
\usepackage{stmaryrd}
\usepackage[numbers]{natbib}
\usepackage{fullpage}   
\usepackage[linesnumbered,ruled]{algorithm2e}

\usepackage[foot]{amsaddr}

\numberwithin{equation}{section}
\numberwithin{figure}{section}

\graphicspath{ {images/} }

\newcommand*\diff{\mathop{}\!\mathrm{d}}

\newtheorem{definition}{Definition}[section]
\newtheorem{remark}[definition]{Remark}

\title{Modelling and Analysis of Non-Contacting Mechanical Face Seals with Axial 
Disturbances and Misalignment}

\author{Ben Ashby\textsuperscript{1}}

\author{Tristan Pryer\textsuperscript{1,2}}

\author{Nicola Bailey\textsuperscript{3,*}}

\address[1]{Institute for Mathematical Innovation, University of Bath, Bath, UK.}
\address[2]{Department of Mathematical Sciences, University of Bath, Bath, UK }
\address[3]{Department of Engineering, King's College London, London, UK.}
\address[*]{Corresponding Author, email nicola.bailey@kcl.ac.uk}

\begin{document}

\begin{abstract}
  Advancements in industrial applications are driving developments in non-contacting
  mechanical seal technology. Key requirements include improvements in efficiency and
  reliability, which lead to smaller clearances, lower frictional losses, and
  minimisation of wear, during operation. However, a critical consideration is the
  effect of external disturbances experience by the seal from the local environment
  which may cause destabilisation, and lead to premature failures through unanticipated
  face contact. 
  
  This work examines the dynamic behaviour of a non-contacting mechanical face seals,
  where a thin fluid film separates a pair of coaxial discs; a rotor (rotating face) and
  stator (stationary face). It is assumed that the rotor-stator have an angular
  misalignment and operation is under conditions involving large axial disturbances,
  representing external disturbances. A fully coupled unsteady mathematical
  representation is developed, where the fluid flow is coupled to the structural
  response of the stator, and the rotor motion is prescribed. The stator is modelled as
  a spring-mass-damper system, and the fluid film model is based on a lubrication
  approximation of the Navier-Stokes equations. The governing equations are solved via a
  numerical technique based on finite element and Runge-Kutta methods. 
  
  A parameter study reveals the impact of misalignment on the seal dynamics, when
  experiencing an external disturbance. The angle of misalignment and corresponding
  amplitude of forcing can be identified when the minimum fluid film thickness becomes
  less than a given tolerance. This provides insights into safe operating conditions and
  manufacturing tolerances, with the research aiding to improve the design critera and
  reliability of non-contacting mechanical face seals.
\end{abstract}

\maketitle

\section{Introduction}
\label{section:intro}

Traditional mechanical face seals are employed in rotating machinery to minimise
undesired fluid flow (leakage) between locations experiencing differential pressures,
whilst allowing power transmission across the seal interface via a rotating shaft.
Current designs consist of a pair of co-axial discs, where one rotates and the other is
stationary, and maintain continuous contact to minimise leakage. However, this
continuous contact results in wear and a reduced performance over time, leading to a
limited lifespan. Consequently, there are increased costs associated with regular
maintenance and replacements together with increased machine downtime.

Non-contacting mechanical face seals could mitigate these issues, offering improved
performance, especially in applications requiring low frictional losses. These seals
employ a thin fluid film between the rotating face (rotor) and stationary face (stator)
to maintain a clearance. The design aims to capitalise on the local film dynamics,
providing sufficiently stiff and responsive fluid-lubricated capabilities to maintain an
acceptable clearance between the rotor and stator. However, non-contacting face seal
technology is often described as unpredictable and unreliable due to a high number of
premature and unexpected failures. This behaviour is likely attributed to unanticipated
face contact, which occurs when the fluid film fails to sustain the face dynamics.
Underlying causes may include minor manufacturing imperfections in the seal components,
relatively large misalignments between the faces, or external disturbances, which in
turn produce poor rotor-stator tracking and rotor dynamic instabilities
\cite{varney2017impact}. Therefore, understanding the seal behavior under these
undesirable operating conditions is essential to minimise potential seal failures.

Modelling the complex behavior of a non-contacting mechanical face seal requires
coupling the fluid film dynamics to the structural behaviour. Previous models have been
developed to predict the behaviour of the face seal. Etsion demonstrated the importance
of the squeeze effect in face seals, where approximations of the fluid film
characteristics are computed analytically from the Reynolds equation, derived using the
narrow seal approximation \cite{etsion1980squeeze}. By coupling the fluid film pressure
to the equation of motion of a flexibly mounted stator seal, Etsion studied the dynamics
of a coned face seal \cite{etsion1982dynamic}. The stability of a coned bearing with
similar geometry and incompressible flow was investigated by Green and Etsion, who
provided a closed-form nonlinear numerical solution for the bearing dynamics
\cite{green1983fluid}. Subsequently, Green and Barnsby formulated a compressible
Reynolds equation model for the fluid film coupled to the structure, arranged into a
single state space form, allowing the fluid film lubrication and the dynamics to be
solved concurrently \cite{green2001simultaneous}. Due to nonlinear effect in the fluid
film, the load carrying capacity has been shown to be greater when transient-state
theory is applied, compared to steady-state \cite{chen2022experimental}. As the
transient-state model is more accurate in predicting experimental seal behaviour
results, it is critical these effects are considered in modelling. The dynamic stability
of dry gas seals, with regards to the structure and operating conditions have been
examined, which highlights the importance of the gas film stiffness and damping
\cite{zhang2024research}.

More recently, the analytical formulation of the coned-face seal was enhanced by
combining the elastic-plastic contact model during steady-state operation to support the
hypothesis that premature seal failure is attributed to undesired face contact
\cite{varney2017impact}. A flexibly mounted rotor geometry was examined by Green, where
the analytical expressions for the stiffness and damping coefficients of the fluid film
were given based on the small perturbation assumption \cite{Green1986rotor}. These
studies provide an overview of the general behaviour under standard operating
conditions. Modified surfaces on the seal faces can be used in an effort to improve the
fluid film characteristics. Studies have examined the effect a conical face, spiral
grooved face, wave face and radial grooved face on the pressure distribution in the
fluid film \cite{blasiak2016parametric, chen2024prediction}. 

However, during operation, the seal may experience non-ideal operating conditions, such
as external disturbances. Garratt et al. examined the dynamics of the coupled
fluid-structure interaction of a high-speed air-lubricated bearing (compressible flow)
with parallel faces \cite{garratt2010compressible}. The rotor had prescribed periodic
axial oscillations with an amplitude smaller than the equilibrium fluid film thickness,
and the stator could move axially in response to the induced film dynamics. Following a
similar approach, Bailey et al. considered the case of a fluid-lubricated bearing with
incompressible flow \cite{bailey2014dynamics} and compressible flow
\cite{bailey2017dynamics} with an axisymmetric coned rotor shape and examined the
extreme case of rotor axial disturbances with larger amplitudes than the equilibrium
film thickness. Results indicate that in the case of parallel faces, the lubrication
force prevents contact; however, the fluid gap can become very small, potentially
invalidating the classical no-slip velocity condition. Extending the analysis to
incorporate a slip boundary condition of the Navier type indicates that parallel
bearings do not experience face contact even though the bearing gap can become very
small, whereas a coned bearing can have possible face contact for certain critical
values of the magnitude of the rotor oscillation, conical angle, and slip condition
\cite{bailey2015dynamics,bailey2016dynamics}. 

In practice, when the seals are assembled, the faces may become misaligned, resulting in
a relative tilt angle and loss of axisymmetry. The lubrication regime transition during
the startup period of a non-contacting mechanical seal has been studied experimentally
with small axial vibrations, and angular misalignment showing serious leakage can occur
during unpredicted contact \cite{xu2022transient}. Therefore, the misalignment can
significantly affect the seal dynamics, especially when exposed to large external
disturbances. To model the non-axisymmetric case, the Reynolds equation must be
approximated numerically, whereas exact solution could be derived for the axisymmetric
geometry. This leads to a significant increase in computational cost relative to
axisymmetric problems. However, approximating the pressure numerically allows
significant flexibility in the sense that different bearing geometries and potentially
different fluid models could be accommodated with little difficulty. The dynamic
behavior of other seal types has been examined for non-axisymmetric geometry under
standard operating conditions using various numerical techniques, including the finite
element, finite volume, and finite difference methods (FEM, FVM, and FDM). For example,
the FEM and FVM were utilized for the dynamic analysis of spiral groove gas face seals
\cite{miller2001numerical}, the FVM was employed by Blasiak to examine the impulse gas
face seal with grooves on both faces \cite{blasiak2019numerical}, and the FDM was
utilised by Li et al. for studying spiral groove liquid film seals \cite{li2020dynamic}.

This work focuses on a non-contacting mechanical face seal that incorporates
imperfections in the form of misalignment of the rotor and stator, leading to a
non-axisymmetric fixed tilted rotor geometry. It is assumed that the surfaces are free
from artifacts. The seal will be examined under extreme operating conditions, namely
large axial rotor disturbances. The outline of the rest of the article is as follows.
The mathematical model for the fully coupled unsteady system is derived in Section
\ref{sec:mathmodel}. The fluid film model is based on thin film flow (lubrication
approximation) of the Navier-Stokes equations, which can be solved to provide analytical
expressions for the fluid velocity and a modified Reynolds equation. The rotor has
prescribed axial oscillations, whereas the stator is free to move axially, modeled as a
spring-mass-damper system, in response to the fluid film dynamics. The stator and fluid
film are coupled through the force from the fluid film on the stator. The numerical
model is formulated in Section \ref{sec:numtech} based on the finite element method for
the modified Reynolds equation and a fourth-order Runge-Kutta method for the stator
equation. A parameter study is carried out in Section \ref{sec:results} to identify
possible destabilising behaviour and investigate the effect of face misalignment in a
non-contacting mechanical face seal. This includes the case of large amplitude
oscillations of the rotor used to simulate potential destabilising behavior of extreme
operating systems. The goal is that by studying the effect of face misalignment, the
resulting system characteristics could aid designers in identifying safe operating
conditions, increasing seal efficiency, and identifying suitable manufacturing
tolerances.

\section{Mathematical model}
\label{sec:mathmodel}

\begin{figure}[h!]
\begin{center}
\def\svgwidth{160mm}
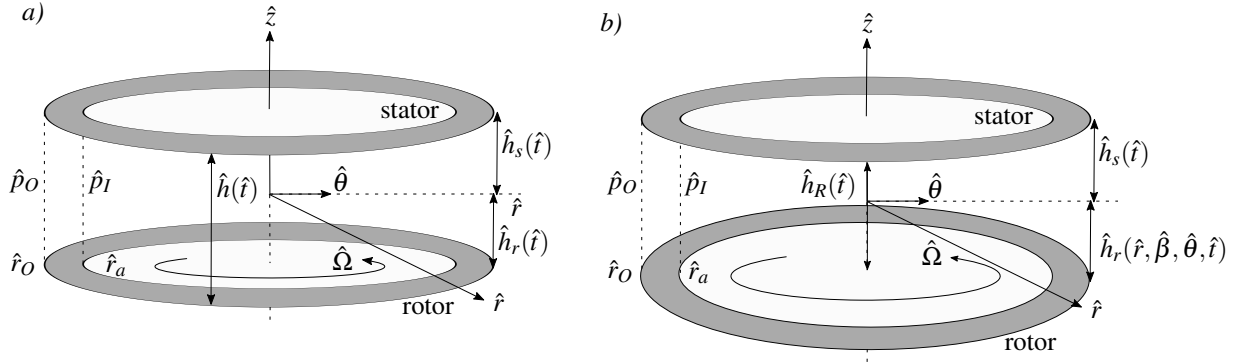 
\end{center}
\caption{Geometry of a fluid-lubricated thrust bearing for a) parallel rotor and stator
and b) tilted rotor in a cylindrical polar coordinate system with tilt angle
$\hat{\beta}$ in cylindrical coordinate system (${\hat{r},\hat{\theta},\hat{z}}$).}
\label{fig:geometry}
\end{figure}

Industrial bearing designs are complex, therefore a simplified mathematical model is
developed by retaining the key features of a fluid-lubricated thrust bearing geometry as
shown in Figure \ref{fig:geometry} for a parallel and titled rotor configuration. The
rotor and stator are modelled as a pair of coaxial annuli with smooth faces, separated
by a thin fluid film. For typical applications the rotor is mounted to a high speed
shaft, which has relative angular rotation $\hat{\Omega}$, whereas the stator is
attached via springs to a housing which is stationary relative to the shaft and 
rotor. 		

In practice the pressure at the inner and outer radius is due to conditions from the
wider system that the bearing is contained within; pressure $\hat{p}_I$ is imposed at
the inner radius ($r_I$) and $\hat{p}_O$ at the outer radius ($r_O$), allowing a
differential pressure to drive a radial flow. A non-axisymmetric configuration is
considered where the rotor and stator can move axially with relative positions
$\hat{z}=\hat{h}_r$ and $\hat{z}=\hat{h}_s$, respectively, in a cylindrical polar
coordinate system $(\hat{r},\hat{\theta},\hat{z})$. The stator height is given by
$\hat{h}_s(t)$ and is assumed to be perpendicular to the $z$ axis and has a tilt angle
$\hat{\beta}$ which is fixed in time and is assumed to be small such that
$\sin\hat{\beta}\sim\hat{\beta}$ and $\cos\hat{\beta}\sim1$. The rotor height is defined
by 
\begin{equation}\label{eqn:dimhr}	
\hat{h}_r(\hat{r},\hat{\beta},\hat{\theta},\hat{t})
=
\hat{h}_R(\hat{t}) + \hat{r}\tan\hat{\beta}\sin\hat{\theta}\,,
\end{equation}
where $\hat{h}_R(\hat{t})$ denotes the rotor height at $\hat{r}=0$, that is, the height
(i.e. $\hat z$-coordinate) at the centre of the rotor and the rotor pivots to give a
tilt angle $\beta$ about the axis $\hat \theta=0$ and $\hat \theta =\pi$ (i.e. the $\hat
x$-axis of a Cartesian coordinate system). To investigate the dynamics of a
fluid-lubricated thrust bearing in a non-ideal operating environment, the rotor has
prescribed periodic axial oscillations representing vibrations or disturbances which
could act to destabilise the bearing operation. Associated with the oscillations is a
frequency $\hat{f}$ with angular frequency $\hat{\omega}=2\pi \hat{f}$ and an amplitude
$\epsilon \hat{h}_0$, where $\hat{h}_0$ is a typical equilibrium rotor-stator clearance
at the reference height. The time dependent axial position of the rotor is given by
\begin{equation}
\hat{h}_R(\hat{t})
=
\epsilon\hat{h}_0\sin(\hat{\omega}\hat{t})\,.
\label{eqn:dimhR}	
\end{equation}
The corresponding rotor-stator clearance is given by
\begin{equation}
\label{eqn:dimhath}	
\hat{h}(\hat{r},\hat{\beta},\hat{\theta},\hat{t}) 
= 
\hat{h}_s(\hat{t}) - \hat{h}_r(\hat{r},\hat{\beta},\hat{\theta},\hat{t})\,.
\end{equation}
The temporal and spatial dependence in the rotor height are independent of each other.

The classical continuum representation is assumed applicable for modelling the fluid
film which separates the rotor and stator, since a typical equilibrium rotor-stator
clearance is $\hat{h}_0=5\text{x}10^{-5}$ m, several orders of magnitude larger than the
mean free path of air, $\hat{l}= 6.8\text{x}10^{-8}$ m. The fluid is modelled as a
Newtonian viscous fluid characterised by the scalar pressure $\hat{p}$ at time $\hat{t}$
and density field $\hat{\rho}$, where the velocity field associated with the fluid flow
in the given coordinate system is $\hat{\mathbf{u}} = (\hat{u},\hat{v}, \hat{w})$. The
dynamic viscosity $\hat{\mu}$ is assumed constant throughout the flow.

\subsection{Fluid Flow Model}
\label{subsec:airflowmodel}
     
A model for the incompressible fluid flow through the bearing is derived from the
classical Navier Stokes equations. On coupling the rotor-stator clearance to the
pressure field for a thin incompressible lubricating film, a modified Reynolds equation
is derived to include a tilted rotor geometry.

The governing Navier-Stokes momentum equations for incompressible flow are given by
\begin{equation}\label{eqn:NS}
    \hat{\rho}
    \left(
         \frac{\partial \hat{\bf{u}}}{\partial \hat{t}}
         +
         (\hat{\bf{u}}\cdot\hat{\nabla})\hat{\bf{u}}
    \right) 
    =
    -\hat{\nabla} \hat{p} 
    +
    \hat{\mu}\hat{\nabla}^2\hat{\bf u}
    +
    \hat{\bf{b}}\,,
\end{equation}
where $\hat{\mathbf{b}}$ is any external body forces per unit, including gravity
$\hat{\mathbf{g}}$. Full details of the derivation are given in Batchelor
\cite{batchelor2000introduction}. This must be coupled to the continuity equation for
incompressible flow, given by 
\begin{equation}\label{eqn:cont}
    \hat{\nabla}\cdot {\bf{\hat{u}}}=0\,.
\end{equation}

\subsubsection{Boundary Conditions}
\label{subsubsec:velbc}

Taking the stator and rotor as solid surfaces, a no-slip and no-flux velocity boundary
condition must be imposed on both these surfaces. The rotor has a constant azimuthal
velocity and is considered displaced in the axial direction with prescribed
oscillations, whereas the stator is able to move axially in response to the interaction
with the fluid. This results in velocity boundary conditions  
\begin{equation}\label{eqn:dimvelbc}
  \begin{alignedat}{4}
      \hat u &= 0, \qquad && 
      \hat v = \hat \Omega \hat r, \qquad && 
      \hat w = \frac{\partial \hat{h}_r}{\partial \hat t},  \qquad & 
      \text{  at  } & \hat z = \hat{h}_r, \\
      \hat u &= 0, \qquad && 
      \hat v = 0, \qquad && 
      \hat w = \frac{\diff \hat{h}_s}{\diff \hat t}, \qquad & 
      \text{  at  } & \hat z = \hat{h}_s,
  \end{alignedat}
\end{equation}
for the rotor and stator, respectively, if the rotor and stator were perfectly aligned
(that is, $\hat{\beta} = 0$). To find the correct boundary condition, the rotational
component of the above velocity field should be rotated about the $\hat{x}$-axis (the
axial movement of the rotor is assumed to remain parallel to the $\hat{z}$-axis). This
gives the misaligned rotor boundary conditions as

\begin{equation}\label{eq:bc_rotor_dimensional}
    \begin{alignedat}{4}
      \hat{u} 
      &= - \hat{\Omega} \hat{z} \hat{\beta}\cos \hat{\theta}, 
      &\quad
      \hat{v} 
      &= \hat{\Omega} \hat{r} + \hat{\Omega} \hat{z} \hat{\beta} \sin \hat{\theta}, 
      &\quad
      \hat{w} 
      &= \hat{\Omega} \hat{r}  \hat{\beta} \cos \hat{\theta} 
      + 
      \frac{\partial \hat{h}_r}{\partial \hat{t}},
      &\quad
      \text{  at  } & \hat z = \hat{h}_r, \\
      \hat{u} 
      &= 0,
      &
      \hat{v} 
      &= 0, 
      &
      \hat{w} 
      &= \frac{\diff \hat{h}_s}{\diff \hat t},
      &
      \text{  at  } & \hat z = \hat{h}_s,
    \end{alignedat}
\end{equation}
after using the small angle approximations $\sin \hat{\beta} \approx \hat{\beta}$ and
$\cos \hat{\beta} \approx 1$ (full details are given in \autoref{sec:bcs_stuff}).
		
A differential pressure is imposed across the bearing to drive a radial flow. Pressure
boundary conditions are defined at the inner and outer radii of the bearing as
\begin{equation}
\hat{p}=\hat{p}_I \quad \text{at} \quad  \hat{r}=\hat{r}_I \quad
\text{and} \quad
\hat{p}=\hat{p}_O 	\quad \text{at} \quad \hat{r}=\hat{r}_O\,,
\end{equation}
respectively.

\subsubsection{Dimensional Analysis}
\label{subsubsec:dim}

To determine the governing physical processes at leading order, the Navier-Stokes
equations (\ref{eqn:NS})-(\ref{eqn:cont}) are expressed in dimensionless variables. A
typical bearing pressure and radius is taken to be $\hat{P}$ and $\hat{r}_0$,
respectively, with dimensionless time variable $t=\hat{t}/\hat{T}_t$. The dimensionless
pressure is defined by $p=\hat{p}/\hat{p}_a$ and dimensionless radial, azimuthal and
axial variables are $r=\hat{r}/\hat{r}_0$, $\hat{\theta}=\theta$ and
$z=\hat{z}/\hat{h}_0$, respectively. The non-dimensional velocity components are given
by $u=\hat{u}/\hat{U}$, $v=\hat{v}/\hat{V}$ and $w=\hat{w}/\hat{W}$ for typical radial,
azimuthal and axial velocities $\hat{U}$, $\hat{V}=\hat{\Omega} \hat{r}_0$ and
$\hat{W}=\hat{h}_0/\hat{T}_t$, respectively. The tilt angle is assumed to be small
giving the non-dimensional angle $\beta=\hat{\beta}/\delta_0$.
			
Non-dimensional groups include the radial and azimuthal Reynolds numbers, the Reynolds
number ratio and a corresponding scaling given by
\begin{equation}\label{eqn:reynoldsnumber}
Re_U = \frac{\hat{h}_0\hat{U}}{\hat{\nu}}, \quad 
Re_\Omega = \frac{\hat{r}_0^2\hat{\Omega}}{\hat{\nu}}, \quad
Re^r = \frac{Re_\Omega}{Re_U} = \frac{\hat{\Omega} \hat{r}_0}{\hat{U}\delta_0}\quad
\text{and}\quad
Re^* = Re^r\delta_0 = \frac{\hat{\Omega} \hat{r}_0}{\hat{U}},
\end{equation}
respectively, characterising the dimensionless flow properties. The aspect ratio
$\delta_0$, squeeze number $\tilde{\sigma}$ and Froude number $Fr$ are defined as		
\begin{equation}
\delta_0 
= 
\frac{\hat{h}_0}{\hat{r}_0},
\quad
\sigma
=
\frac{\hat{r}_0}{\hat{U}\hat{T}}
\quad
\text{and}
\quad
Fr
=
\frac{\hat{U}}{(\hat{g}\hat{h}_0)^{-\frac{1}{2}}},
\label{eqn:deltasigmaFr}
\end{equation}
respectively, where $\hat{g}$ is the acceleration due to gravity.
			
For thin film bearings the aspect ratio is small, $\delta_0 \ll 1$ and effects of
inertia are neglected due to the reduced Reynolds number $Re_U \delta_0^2\ll 1$, with
terms of the order $Re_U\delta_0 (Re^*)^2$ are considered to be negligible, with
$(Re^*)^2$ of $O(1)$ at most. The squeeze number $\tilde{\sigma}$, taken to be of
$O(1)$, characterises any time dependent effects whilst the Froude number $Fr$
parametrises the importance of the gravitational effects relative to the radial flow.
However gravity can be neglected with $Re_U \delta_0Fr^{-2}\ll 1$; this is consistent
with lubrication theory provided the Froude number is $O(1)$. The effects of viscosity
are retained at leading order with the pressure scaled as
$P=\mu\hat{r}_0\hat{U}/{\hat{h}_0}^2$.

Applying the above scalings to the Navier-Stokes momentum equations (\ref{eqn:NS}), the
radial component becomes 
\begin{subequations}
\begin{equation}
    \begin{split}
Re_U\delta_0&
\left(
    \sigma\left(\frac{\partial u}{\partial t}+ w\frac{\partial u}{\partial z}\right) 
    + 
    u\frac{\partial u}{\partial r}+Re^*\frac{v}{r}\frac{\partial u}{\partial \theta} 
    -
    \left(Re^*\right)^2\frac{v^2}{r}
\right)
=\\
- &\frac{\partial p}{\partial r}+{\delta_0}^2
\left(
    \frac{1}{r}\frac{\partial}{\partial r}\left(r\frac{\partial u}{\partial r}\right) 
    + 
    \frac{1}{r^2}\frac{\partial^2 u}{\partial \theta^2} - \frac{u}{r^2}
\right) 
+ 
\frac{\partial^2 u}{\partial z^2} 
- 
2 Re^*{\delta_0}^2\frac{1}{r^2}\frac{\partial v}{\partial \theta},
\label{eqn:NDNScylu}
    \end{split}
\end{equation}
the azimuthal component becomes
\begin{equation}
    \begin{split}
Re_U\delta_0Re^*&
\left(
    \sigma\left(\frac{\partial v}{\partial t}+ w\frac{\partial v}{\partial z}\right) 
    + 
    u\frac{\partial v}{\partial r}  
    + 
    \frac{uv}{r}
\right) 
+ 
Re_U\delta_0(Re^*)^2\frac{v}{r}\frac{\partial v}{\partial \theta}
=\\
-&\frac{1}{r}\frac{\partial p}{\partial \theta}+Re^*{\delta_0}^2
\left(
    \frac{1}{r}\frac{\partial}{\partial r}\left(r\frac{\partial v}{\partial r}\right) 
    + 
    \frac{1}{r^2}\frac{\partial^2 v}{\partial \theta^2} + \frac{v}{r^2}
\right) 
+ 
Re^*\frac{\partial^2 v}{\partial z^2}
+ 
2{\delta_0}\frac{1}{r^2}\frac{\partial u}{\partial \theta},
\label{eqn:NDNScylv}
    \end{split}
\end{equation}
and the axial component becomes
\begin{equation}
    \begin{split}
Re_U{\delta_0}^3\sigma&
\left(
    \sigma\left(\frac{\partial w}{\partial t}+ w\frac{\partial w}{\partial z}\right) 
    + 
    u\frac{\partial w}{\partial r} 
    + 
    Re^*\frac{v}{r}\frac{\partial w}{\partial \theta}
\right) 
=\\
&Re_U\delta_0{Fr}^{-2} - \frac{\partial p}{\partial z}
+
\sigma{\delta_0}^4
\left(
    \frac{1}{r}\frac{\partial}{\partial r}\left(r\frac{\partial w}{\partial r}\right) 
    + 
    \frac{1}{r^2}\frac{\partial^2 w}{\partial \theta^2} 
    + 
    \frac{1}{{\delta_0}^2} \frac{\partial^2 w}{\partial z^2}
\right).
\label{eqn:NDNScylw}
    \end{split}
\end{equation}
\label{eqn:NDNScyl}
\end{subequations}
Similarly the continuity equation becomes
\begin{equation}
\frac{1}{r}\frac{\partial}{\partial r}(ru)
+
Re^*\frac{1}{r}\frac{\partial v }{\partial \theta} 
+ 
\sigma\frac{\partial w}{\partial z}=0\,.
\label{eqn:NDNScylcont}
\end{equation}			
						
Using the above scalings and dimensionless parameters the leading order Navier-Stokes
momentum equations (\ref{eqn:NDNScyl}), where terms $O(\delta_0)$ are neglected, become
\begin{equation}\label{eq:leading_order_momentum}
    -\frac{\partial p}{\partial r}+\frac{\partial^2 u}{\partial z^2}=0,
    \quad
    -\frac{1}{r}\frac{\partial p}{\partial \theta} 
    + 
    Re^* \frac{\partial^2 v}{\partial z^2}
    =
    0,
    \quad
    \frac{\partial p}{\partial z}=0.
\end{equation}
Similarly the continuity equation (\ref{eqn:NDNScylcont}) becomes 
\begin{eqnarray}
    \frac{1}{r}\frac{\partial}{\partial r}(r u) 
    + 
    Re^*\frac{1}{r}\frac{\partial v}{\partial \theta} 
    + 
    \sigma\frac{\partial w}{\partial z}=0,
    \label{eqn:momcont}
\end{eqnarray}
using the definition of the squeeze number in (\ref{eqn:deltasigmaFr}). 
			
The velocity boundary conditions in (\ref{eqn:dimvelbc}) are given by	

\begin{equation}\label{eq:noslipnondimvelbc}
  \begin{alignedat}{4}
          u &= 0 \qquad && v = \Omega r \qquad && w = \frac{\partial h_r}{\partial t} 
          + \frac{Re^*}{\sigma} r \beta \cos \theta \qquad & \text{  at  } & z = h_r, \\
          u &= 0 \qquad && v = 0 \qquad && 
          w = \frac{\diff h_s}{\diff t} \qquad & \text{  at  } & z = h_s,
  \end{alignedat}
\end{equation}
in dimensionless variables for the rotor and stator, respectively (further details are
provided in \autoref{sec:appendix_dimensional}).

In dimensionless variables the pressure boundary conditions become
\begin{equation}\label{eqn:pressurebc}
p=p_I \quad \text{at} \quad r=a, \quad \text{and}\quad p=p_O \quad \text{at} \quad r=1, 
\end{equation}
with dimensionless pressures $p_I=\hat{p}_I/\hat{P}$ and $p_O = \hat{p}_O/\hat{P}$. 

The rotor height together with the rotor clearance from equations (\ref{eqn:dimhr}) 
are 	
\begin{equation}\label{eqn:dimhrtheta}		
    h_r(r,\beta,\theta,t)=\epsilon\sin t+ r\beta\sin\theta,
\end{equation}
taking a time scale of $\hat{T}=1/{\omega}$.
The corresponding rotor-stator clearance is given by
\begin{equation}\label{eqn:dimh}	
    h(r,\beta,\theta,t)=h_s(t) - h_r(r,\beta,\theta,t).
\end{equation}

\subsection{Modified Reynolds Equation}
\label{subsec:Reyneqn}

Analytical responses for the velocity components $(u,v,w)$ can be found in terms of the
unknown stator position $h_s(t)$ by integrating the governing equations
(\ref{eq:leading_order_momentum})-(\ref{eqn:momcont}) and applying the dimensionless velocity boundary
conditions (\ref{eq:noslipnondimvelbc}). The radial, azimuthal and axial velocities are
given, respectively, by
\begin{subequations}
\begin{eqnarray}
    u(r,\beta,\theta,z,t) 
    &=& 
    \frac{1}{2}\frac{\partial p}{\partial r}\left(z-h_s\right)\left(z-h_r\right),
    \label{eqn:uvel}\\			
    v(r,\beta,\theta,z,t)
    &=&
    \frac{1}{2Re^*r}\frac{\partial p}{\partial\theta}\left(z-h_s\right)
    \left(z-h_r\right) - \frac{r}{h}\left(z-h_s\right),\label{eqn:vvel}\\
    w(r,\beta,\theta,z,t) 
    &=&
    \frac{\partial h_r}{\partial t} + \frac{\partial h_r}{\partial \theta} 
    - 
    \frac{1}{\sigma r}\frac{\partial}{\partial r}
    \left(\frac{r}{12}\frac{\partial p}{\partial r}(z-h_r)^2(2z+h_r-3h_s)\right)
    \label{eqn:wvel}\\
    &&
    - 
    \frac{1}{\sigma r}\frac{\partial}{\partial \theta}
    \left(\frac{r}{12}\frac{\partial p}{\partial \theta}(z-h_r)^2(2z+h_r-3h_s) 
    - \frac{r}{2h}(z-h_r)(z+h_r-h_s)\right).\nonumber
\end{eqnarray}
\end{subequations}
The tilt angle of the rotor and $\theta$ dependence is implicitly expressed in the rotor
height $h_r$, rotor-stator clearance $h$ and pressure $p$.
			
To derive the modified Reynolds equation, the continuity equation (\ref{eqn:momcont}) is
integrated between the rotor and stator, and the Leibniz integral rule for
differentiation under the integral sign is applied to give

\begin{equation}\label{eqn:noslipslipreysame}
  \begin{split}
    w(h_s)-w(h_r)
    &=
    -\frac{1}{\sigma r} 
    \left(\int_{h_r(r)}^{h_s}\frac{\partial}{\partial r}(ru) \diff z 
    + 
    Re^* \int_{h_r(\theta)}^{h_s}\frac{\partial v}{\partial \theta}\diff z\right) \\
    &=
    -\frac{1}{\sigma r} 
    \left(
    \frac{\partial}{\partial r}\int_{h_r}^{h_s}ru \diff z + 
    \left(\frac{\partial h_r}{\partial r} (ru) \right)_{z=h_r}
    +
    Re^*\frac{\partial}{\partial \theta}\int_{h_r}^{h_s}v \diff z 
    + 
    \left(\frac{\partial h_r}{\partial \theta} v \right)_{z=h_r}
    \right).
  \end{split}
\end{equation}
Applying the velocity boundary conditions in (\ref{eq:noslipnondimvelbc}), gives
\begin{eqnarray}\label{eq:pre_reynolds}
    \sigma\left(\frac{\diff h_s}{\diff t}-\frac{\partial h_r}{\partial t}\right)
    -
    Re^* r \beta \cos \theta
    =
    -
    \frac{1}{r}\frac{\partial}{\partial r}\int^{h_s}_{h_r}(ru)\diff z
    -
    \frac{Re^*}{r}\frac{\partial}{\partial \theta}\int^{h_s}_{h_r}v\diff z,
\end{eqnarray}
and substituting the expressions for the radial and azimuthal velocities
(\ref{eqn:uvel}-\ref{eqn:vvel}) in \eqref{eq:pre_reynolds}, after some manipulations the
modified Reynolds equation is found to be
\begin{equation}\label{eqn:Reynolds}	
    \sigma \frac{\partial h}{\partial t}
    -
    \frac{1}{12r}\frac{\partial}{\partial r}
    \left(rh^3\frac{\partial p}{\partial r}\right) 
    - 
    \frac{1}{12r}\frac{\partial}{\partial \theta}
    \left(\frac{h^3}{r}\frac{\partial p}{\partial\theta}\right) 
    - 
    \frac 1 2 r\beta\sin\theta
    =
    0.
\end{equation}
Alternatively, after defining a scaled squeeze number $\tilde{\sigma} = 12 \sigma$ and
recalling the expressions for the gradient and divergence in plane polar coordinates,
the equation can be given in the compact form
\begin{equation}\label{eq:final_Reynolds}
    \tilde{\sigma} \frac{\partial h}{\partial t} 
    - 
    \nabla \cdot \left(h^3 \nabla p\right)
    =
    6 Re^* r \beta \cos \theta.
\end{equation}
A similar method was used when calculating the axial velocity, as it requires the
application of the Leibniz rule. The modified Reynolds equation (\ref{eqn:Reynolds})
expresses the relationship between the pressure $p$ and the rotor-stator clearance $h$,
giving the flow characteristics when solved. 

\begin{remark}[Different bearing designs]
    In the derivation of the modified Reynolds equation, the specific functional form of
    the rotor-stator gap (i.e. a misaligned parallel bearing design) was only used to
    derive the right hand side. This means that other bearing surfaces such as spiral
    grooves can be implemented in the same framework with relative ease.
\end{remark}

\subsection{Modelling Rotor-Stator Dynamics}
\label{subsec:modelrsdynamics}

The fluid flow is coupled to the stator equation through the pressure at the surface of
the stator exerting an axial force which drives its motion. Thus, the axial displacement
of the stator is modelled using Newton's second law as a spring-mass-damper system,
giving the axial position of the stator as
\begin{equation}
\hat{m}\frac{\diff^2\hat{h}_s}{\diff\hat{t}^2}
+
\hat{D}_a\frac{\diff\hat{h}_s}{\diff\hat{t}}
+
\hat{K}_z(\hat{h}_s-\hat{h}_{0s})
=
\hat{F}\left(\hat{t}\right) - \hat{m}\hat{g},
\label{eqn:dimstator}
\end{equation}
where $\hat{m}$ is the mass of the stator and $\hat{h}_{0s}$ the stator equilibrium
height in the absence of gravity. The bearing quantities $\hat{D}_a$ and $\hat{K}_z$ are
the external damping and effective restoring force coefficients, respectively, applied
to constrain the stator to the equilibrium position. The resultant axial force on the
stator is given by 
\begin{equation}\label{eqn:dimforce}
    \hat{F}\left(\hat{t}\right)
    =
    \int_0^{2\pi}\int^{\hat{r}_O}_{\hat{r}_I}\hat{\sigma}_{zj}\hat{n}_{j} \hat{r}
    \diff \hat{r}\diff\hat{\theta}
    =
    \int_0^{2\pi}\int^{\hat{r}_O}_{\hat{r}_I}\left(\hat{p}-\hat{p}_a-2\hat{\mu}
    \frac{\partial \hat{w}}{\partial \hat{z}}\right)\hat{r}\diff \hat{r}
    \diff\hat{\theta},
\end{equation}
due to examining the normal stresses on the stator, with the normal to the stator given
by $\hat{\mathbf{n}}=(0,0,-1)$; a component from the pressure field and component from
the moving fluid. The atmospheric pressure is denoted by $\hat{p}_a$. The stator
equilibrium position is defined as
\begin{equation}
\hat{h}_0=\hat{h}_{0s}-\frac{\hat{m}\hat{g}}{\hat{K_z}},
\end{equation}
using the reference system in Figure \ref{fig:geometry}.

Applying the non-dimensional scalings in Section \ref{subsubsec:dim} the axial stator
deflection from its equilibrium position is given by
\begin{equation}\label{eqn:stator}
    \frac{\diff^2h_s}{\diff  t^2}+D_a\frac{\diff h_s}{\diff t}+K_z(h_s-1)=\alpha F(t),
\end{equation}
with $h_s(t)$ defined at the origin of reference as shown in Figure \ref{fig:geometry}
and the rotor-stator clearance at $\theta=0$ taken to be $h=1$ to give the stator height
$h_s=1$. In (\ref{eqn:stator}), the leading order dimensionless force $F(t)$ is given 
by 
\begin{equation}\label{eq:leading_order_force}
F(t)
=
\int_0^{2\pi}\int^1_a(p-p_a)r\diff r\diff \theta.
\end{equation}
The dimensionless force coupling parameter defined as $\alpha=\hat{\mu}
\hat{U}/\hat{m}\hat{\omega}^2\delta_0^3$, which is taken as $O(1)$. The dimensionless
damping and effective restoring force parameters are given by
$D_a=\hat{D}_a/\hat{m}\hat{\omega}$ and $K_z=\hat{K}_z/\hat{m}\hat{\omega}^2$,
respectively.

\section{Numerical model}\label{sec:numtech}

In this section, a numerical procedure for approximation of the non-dimensional
misaligned bearing problem is described. Let 
\begin{equation}
    \Omega 
    = 
    \{\mathbf{x} \in \mathbb{R}^2 \mid a < |\mathbf{x}| < 1\},
\end{equation} 
where $0 < a <1$, with boundary denoted by $\partial \Omega$. Let the inner boundary be
defined as
\begin{equation}
\Gamma_I:=
    \{\mathbf{x} \in \partial \Omega : |\mathbf{x}| = a\},
\end{equation}
while the outer boundary 
\begin{equation}
    \Gamma_O:=\{\mathbf{x} \in \partial \Omega : |\mathbf{x}| = 1\}.
\end{equation}
The full misaligned bearing problem is given by 

\begin{subequations}
\begin{equation}\label{eqn:Rynlds}
	\tilde{\sigma} \frac{\partial h}{\partial t}
	-
	\frac{1}{r}\frac{\partial}{\partial r}
	\left(rh^3\frac{\partial p}{\partial r}\right) 
	- 
	\frac{1}{r}\frac{\partial}{\partial \theta}
	\left(\frac{h^3}{r}\frac{\partial p}{\partial\theta}\right) 
	= 
	6 r\sin\theta \tan \beta \quad \text{in}\,\, \Omega,
\end{equation}
\begin{equation}\label{eqn:str_hght}
    \frac{\diff^2 h_s}{\diff t^2}
    +
    D_a\frac{\diff h_s}{\diff t}
    +
    K_z(h_s-1)
    =
    \alpha F\bigg(t, h, \frac{\partial h}{\partial t}\bigg), \quad t \in [0,T],
\end{equation}
\begin{equation}\label{eqn:bndry}
	p = p_I \,\text{   on  }\Gamma_I, \quad 
	p = p_O \,\text{   on  }\Gamma_O,
\end{equation}
\begin{equation}\label{eqn:ntl}
	h_s(0) = h^0_s, \quad \frac{\diff h_s(0)}{\diff t} = 0.
\end{equation}
\end{subequations}

In developing a numerical model to solve system \eqref{eqn:Rynlds}-\eqref{eqn:ntl}, it
is assumed that the fluid is quasistatic in the sense that the fluid response to changes
in the rotor-stator clearance $h$ is instantaneous. This is a feature of the leading
order thin film equations \eqref{eq:leading_order_momentum}, which do not contain any
time derivative terms. Time dependence is only introduced through the fluid velocity
boundary conditions \eqref{eq:noslipnondimvelbc}. In what follows, the various aspects
of solving this coupled system numerically are discussed. A numerical solver for the ODE
problem \eqref{eqn:str_hght} is described in \S\ref{sec:solve_ODE} (see also Appendix
\ref{sec:time_stepper}). This ODE is forced by the pressure, determined by Equation
\eqref{eqn:Rynlds}. It is noted that unlike the axisymmetric case (i.e. $\beta = 0)$,
analytical solutions to the Reynolds equation are not readily obtained. As a result, the
numerical ODE solver must be combined with a numerical approximation of the pressure,
which in this work is provided by an adaptive finite element method, described in
\S\ref{sec:FEM} (see also Appendix \ref{sec:fem_appendix}). In \S\ref{sec:adapt} the
algorithm used for mesh adaptivity based upon a local error indicator is described, with
further detail provided in \ref{sec:implement_adaptivity}.

\subsection{Approximation of pressure response}\label{sec:FEM}

In this section the approximation of the forcing function $F$, is addressed. Numerical
solvers, when applied to the ODE \eqref{eqn:str_hght}, will be required to evaluate an
approximation of the forcing function $F$ many times, and for any state of the bearing
system, which at a given time is fully determined by the values of $h_s$, $h_r$ and
$\diff h \slash \diff t$ due to the quasistatic nature of the fluid. This approximation
will be denoted $\mathcal{F}$. Suppose that in a particular state, the stator and rotor
heights are respectively given by $h_{s,\,*}$ and $h_{r,\,*}$, with the rate of change
of the rotor-stator gap given by $\diff h_* \slash \diff t$. Then the spatially variable
gap between rotor and stator is given by 
\begin{equation}
    h_*(r, \theta) 
    =
    h_{s,\,*} - h_{r,\,*} - r \beta\sin\theta.
\end{equation}
These data are sufficient to determine a pressure response $p_*$ which is the solution
of the Reynolds equation

\begin{equation}\label{eqn:time_discrete_reynolds4}
	-
	\frac 1 r \frac{\partial}{\partial r}
	\left(r h_*^3 \frac{\partial p_*}{\partial r}\right)
    -
	\frac 1 r \frac{\partial}{\partial \theta}
	\left(\frac{h_*^3}{r} \frac{\partial p_*}{\partial \theta}\right)
	=
	6 r \beta \sin \theta  
	- 
	\tilde \sigma \frac{\diff h_*}{\diff t}.
\end{equation}
\autoref{eqn:time_discrete_reynolds4} is approximately solved using a piecewise linear
continuous finite element method to obtain the discrete pressure response $\Phi_*$. To
preserve the clarity of the presentation, we provide further details on the spatial
discretisation in Appendix \ref{sec:fem_appendix}. Following this calculation, $\Phi_*$
can be used to obtain an approximation of the leading order pressure force on the stator
given in \autoref{eq:leading_order_force}:

\begin{equation}\label{eq:approx_F}
	\mathcal{F}\left(h_{s,\,*}, h_{r,\,*}, \frac{\diff h_*}{\diff t}\right)
	:= 
	\int_{\Omega'}(\Phi_* - p_a) \diff \mathbf{ x},
\end{equation}
where $\Omega'$ is an approximation of the annular geometry on which the finite element
mesh is defined (see Appendix \ref{sec:fem_appendix} for details).

\subsection{Simulation of Rotor-Stator dynamics}\label{sec:solve_ODE}

The rotor-stator dynamics are governed by the nonlinear second order ordinary
differential equation \eqref{eqn:str_hght}, with initial conditions \eqref{eqn:ntl}.
After defining $x_1 = h_s$ and $x_2 = \diff x_1 \slash \diff t$, this equation may be
expressed as the following first order system of ordinary differential equations

\begin{equation}\label{eqn:stator1order2}
	\begin{split}
		\frac{\diff x_1}{\diff t} &= x_2, \\
		\frac{\diff x_2}{\diff t} &= \alpha F(t, x_1, x_2) - D_a\, x_2 - K_z(x_1-1).
		\end{split}
\end{equation}
This system of equations is nonlinear, with evaluation of $F$ requiring the solution of
the Reynolds equation which depends upon $h$ and its time derivative. Linearisation of
this problem is a significant challenge, rendering implicit methods computationally very
expensive, though they do generally have the benefit of superior numerical stability. We
mention also the possibility of implicit schemes where the Jacobian of $F$ is
approximated numerically. 

In this work, an explicit time-stepping scheme is implemented to solve the system
\eqref{eqn:stator1order2}. This means that approximations of $F$ can be approximated by
$\mathcal{F}$ using the system state at a previous time step as described in
\S\ref{sec:FEM} and Appendix \ref{sec:fem_appendix}. It is emphasised that due to the
stiff nature of the ODE system, care must be taken when solving with explicit schemes,
but this appears to be a reasonable trade-off when compared with the issues faced by
implicit methods in this case. Periodic solutions are not assumed, as a major goal of
this work is to seek the initial response of the bearing dynamics to disturbances, and
to observe the approach of solutions to a periodic cycle. Thus, the stroboscopic map
solver used in e.g. \cite{bailey2017dynamics} is not available. The fourth order
accurate Runge-Kutta method (RK4) 
is chosen to provide fine temporal resolution which may be required when the gap between
rotor and stator becomes small. Further details on the numerical approximation of the
rotor-stator dynamics are provided in Appendix \ref{sec:time_stepper}.

\subsection{Dynamic mesh adaptivity}\label{sec:adapt}

\begin{figure}
	\begin{subfigure}{0.45\linewidth}
	  \includegraphics[trim={5.8cm 3.8cm 5.8cm 4.2cm},clip,width=\linewidth]
      {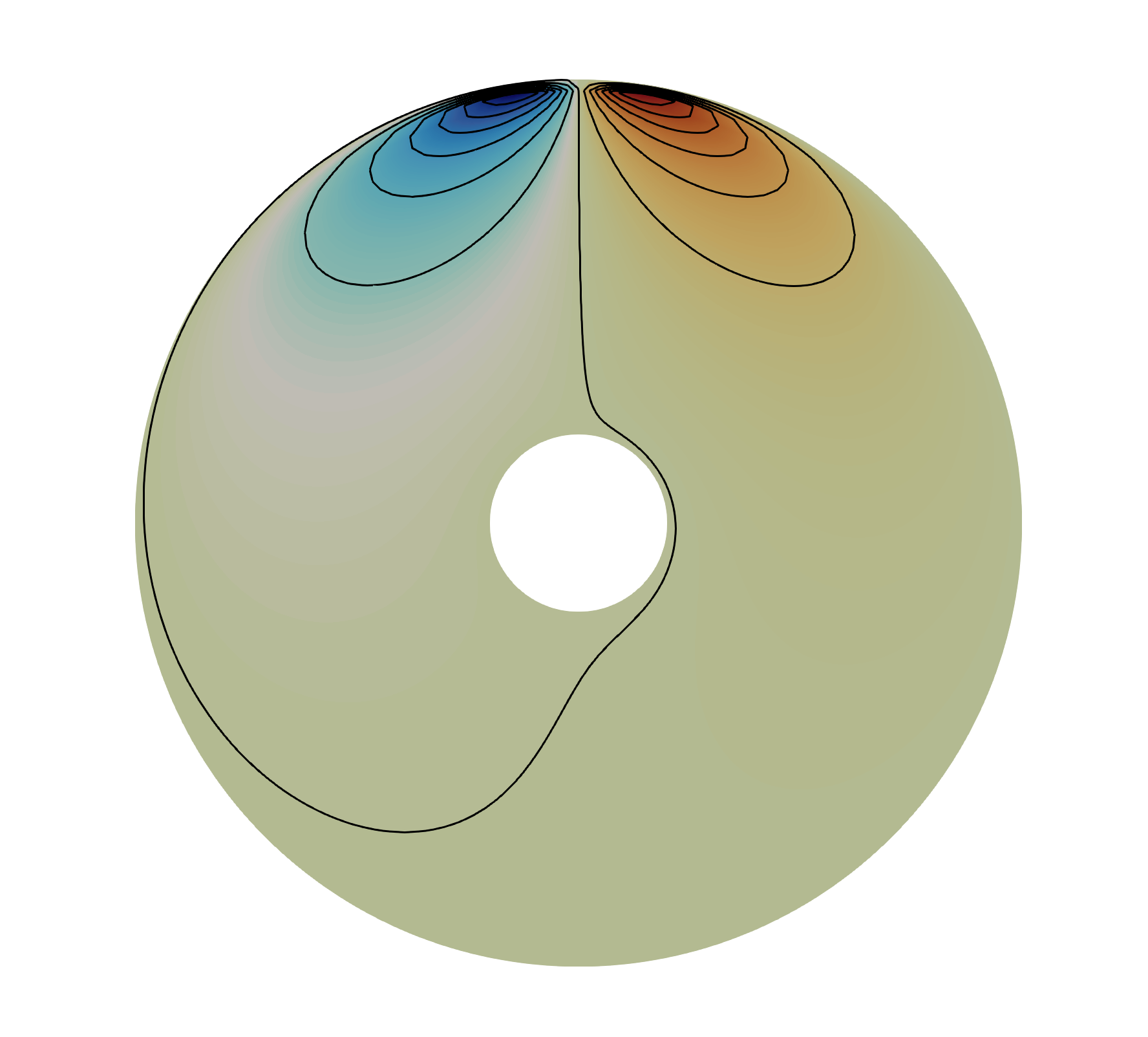}
	  \caption{ Contours of pressure field}
   \label{fig:contours_with_layers}
	\end{subfigure}
	\begin{subfigure}{0.45\linewidth}
		\includegraphics[width=\linewidth, trim={3.3cm 3.3cm 3.3cm 3.3cm},clip]
        {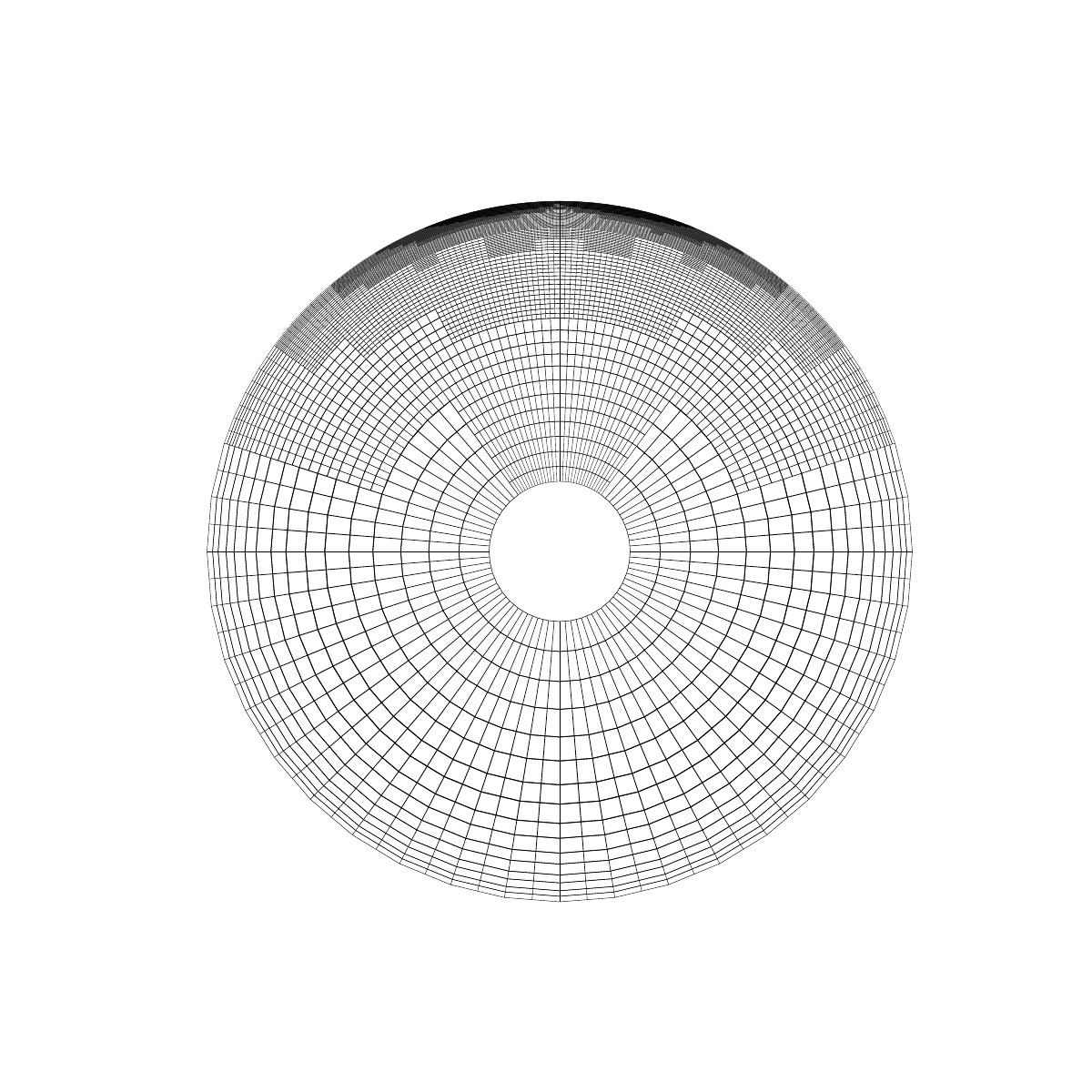}
		\caption{Computational mesh}
  \label{fig:mesh_resolving_layers} 
	  \end{subfigure}
	  \caption{An Illustration of the highly localised pressure peaks that occur when 
      the rotor-stator gap is small, and the adaptive meshing techniques that are used 
      to capture them effectively.}
   \label{fig:adaptive_mesh_and_field}
\end{figure}

Due to the misalignment of rotor and stator, steep pressure layers can occur in the
neighbourhood of the point where the rotor and stator are closest. This behaviour is
observed in Figure \ref{fig:contours_with_layers}, which shows a contour of a pressure
field which results when the stator is disturbed from its equilibrium height by axial
movement of the rotor (see \S\ref{sec:smooth_disturb} for details). A steep pressure
layer between negative (blue) and positive (red) values. This occurs due to the face
clearance narrowing and leads to the Reynolds equation becoming dominated by convection
rather than diffusion (see the discussion in Appendix \ref{sec:pde_stuff}) when the gap
is small. For the layers that arise in such problems, impractically fine meshes are
required to resolve them. One of the major advantages of finite element methods is the
relative ease in which mesh adaptivity can be applied to better capture local features
in the solution, as seen for example in Figure \ref{fig:adaptive_mesh_and_field}. To
address the difficulties described above, dynamic mesh adaptivity was implemented based
upon the well established Kelly error indicator, introduced in
\cite{kelly1983posteriori}. An example of a mesh obtained using this criterion is shown
in Figure \ref{fig:mesh_resolving_layers}, displaying significantly finer mesh
resolution around the layers in the solution, yet keeping the overall number of degrees
of freedom (and hence the computational cost) much lower than would be required on a
uniform mesh sufficient to resolve the solution. Full details on this error indicator
and the implementation of mesh adaptivity are given in Appendix
\ref{sec:implement_adaptivity}.

\begin{figure}
    \centering
    \begin{subfigure}{0.32\linewidth}
    \includegraphics[width=\linewidth]{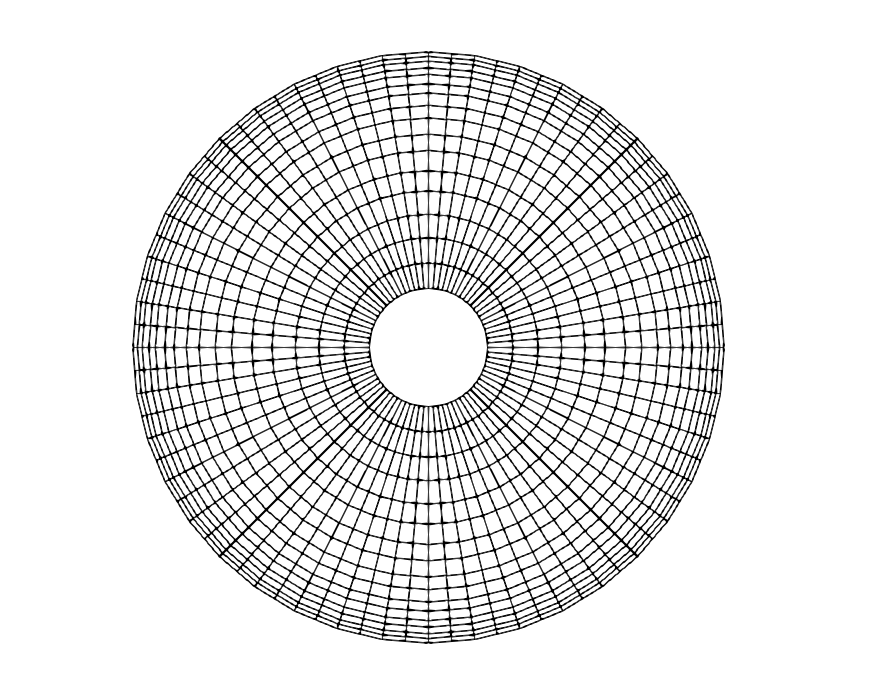}
    \caption{$t = 0.8$}
    \end{subfigure}
    \begin{subfigure}{0.32\linewidth}
    \includegraphics[width=\linewidth]{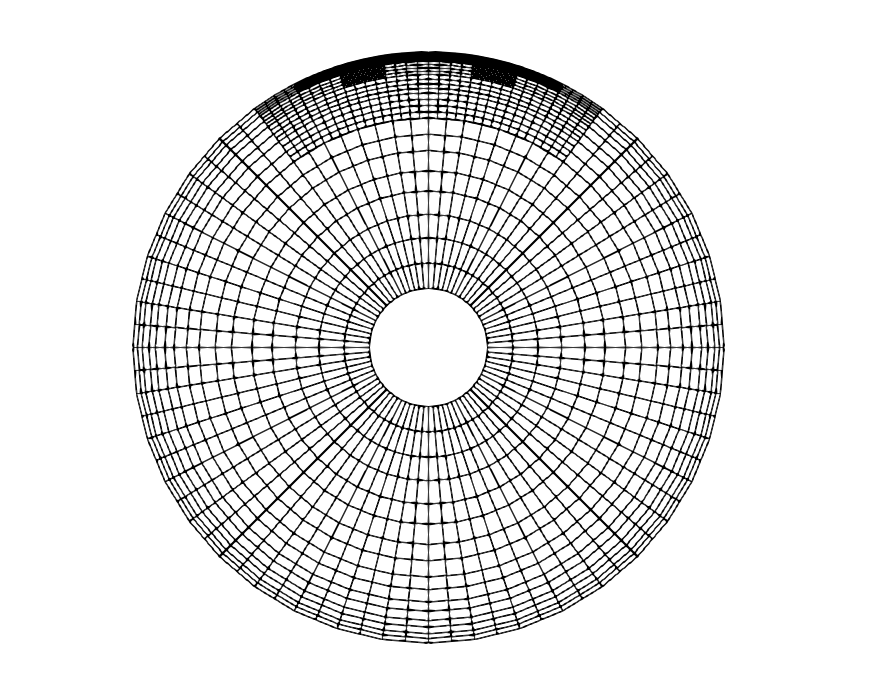}
    \caption{$t = 1.0$}
    \end{subfigure}
    \begin{subfigure}{0.32\linewidth}
    \includegraphics[width=\linewidth]{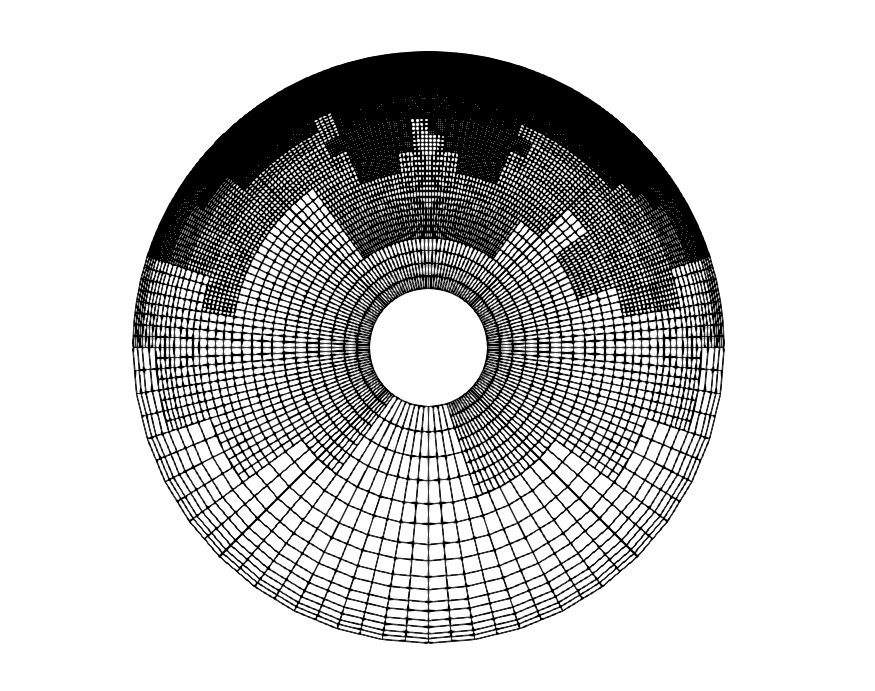} 
    \caption{$t = 1.4$}
    \end{subfigure}\\
    \begin{subfigure}{0.32\linewidth}
    \includegraphics[width=\linewidth]{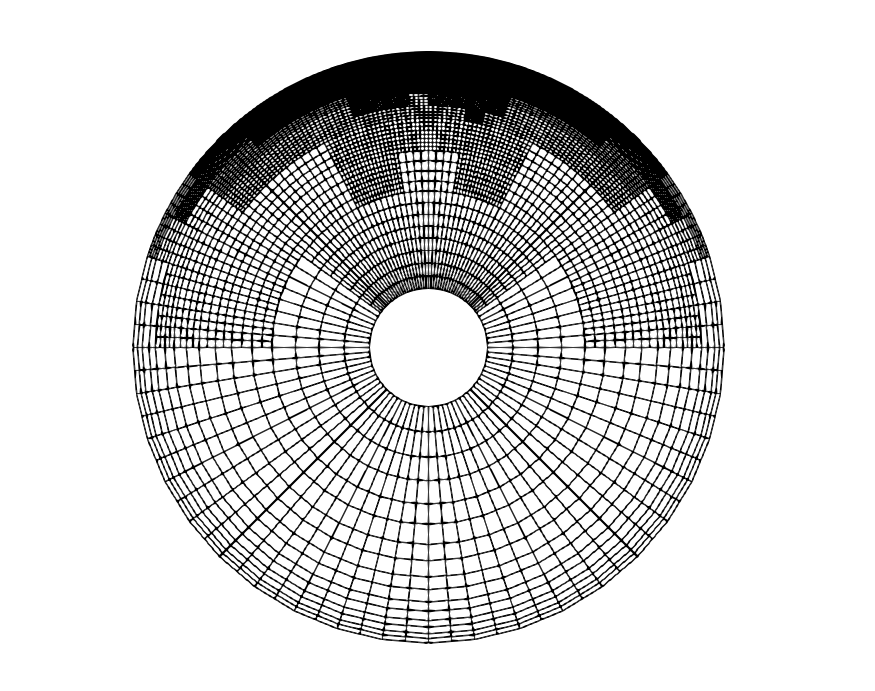}
    \caption{$t = 1.6$}
    \end{subfigure}
    \begin{subfigure}{0.32\linewidth}
    \includegraphics[width=\linewidth]{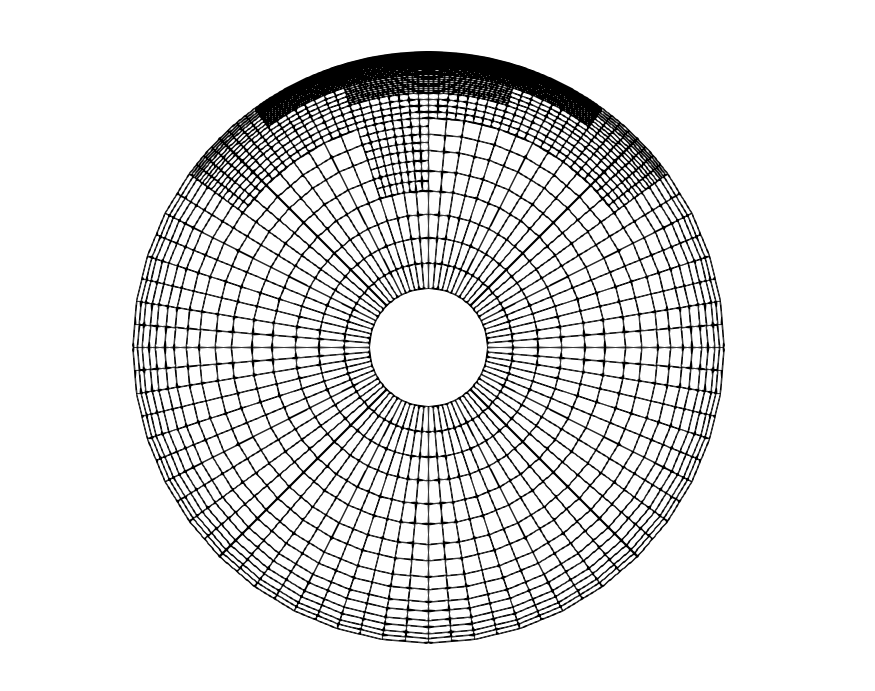}
    \caption{$t = 1.8$}
    \end{subfigure}
    \begin{subfigure}{0.32\linewidth}
    \includegraphics[width=\linewidth]{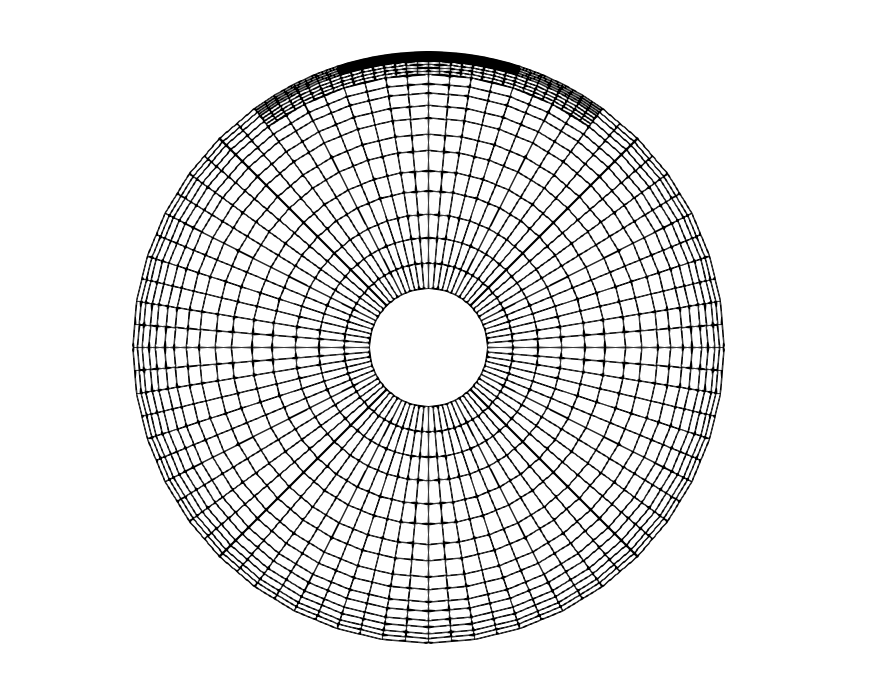}
    \caption{$t = 2$}
    \end{subfigure}
    \caption{A sequence of adaptive meshes, generated automatically during the
    simulation of Example 2, \S\ref{sec:sine_disturb}. The meshes correspond to
    approximately the region of phase $[\tfrac \pi 6, \frac{2 \pi}{6}]$ of the period of
    the prescribed rotor motion during which the rotor and stator are closest. The time
    index of each mesh is given in its respective sub-caption. The mesh is initially
    uniform, and is refined in response to smaller gap $g$ and increasing pressure
    vales. In the bottom row, it is observed that as the region of high pressure becomes
    more and more localised, computational resources are directed to the region where
    the face gap is smallest, leading to a highly locally refined mesh.}
    \label{fig:adaptive_meshing}
\end{figure}

\section{Results}
\label{sec:results}

Numerical results are presented to illustrate the simulated dynamics of a misaligned
thrust bearing under external disturbance. We first examine the relaxation to
equilibrium following a smooth disturbance in \S\ref{sec:smooth_disturb} while in the
following \S\ref{sec:sine_disturb} we model unwanted vibrational disturbances using a
sinusoidal forcing. Finally, in \S\ref{sec:safety} a region of safe operation for the
model thrust bearing is determined by finding critical misalignment angles for various
forcing amplitudes at which contact between rotor and stator occurs. In the following
numerical experiments, Equations \eqref{eqn:stator1order2} are solved with stationary
rotor (i.e. setting $\epsilon = 0$) to first find the equilibrium height of the stator.
This is then used as an initial condition to simulate a disturbance from equilibrium.
Numerical simulations were performed using the \texttt{C++} finite element library
\texttt{deal.II} \cite{dealiiv94}.

\subsection{Example 1: smooth disturbance to equilibrium}
\label{sec:smooth_disturb}

\begin{figure}
    \includegraphics[width=\linewidth]{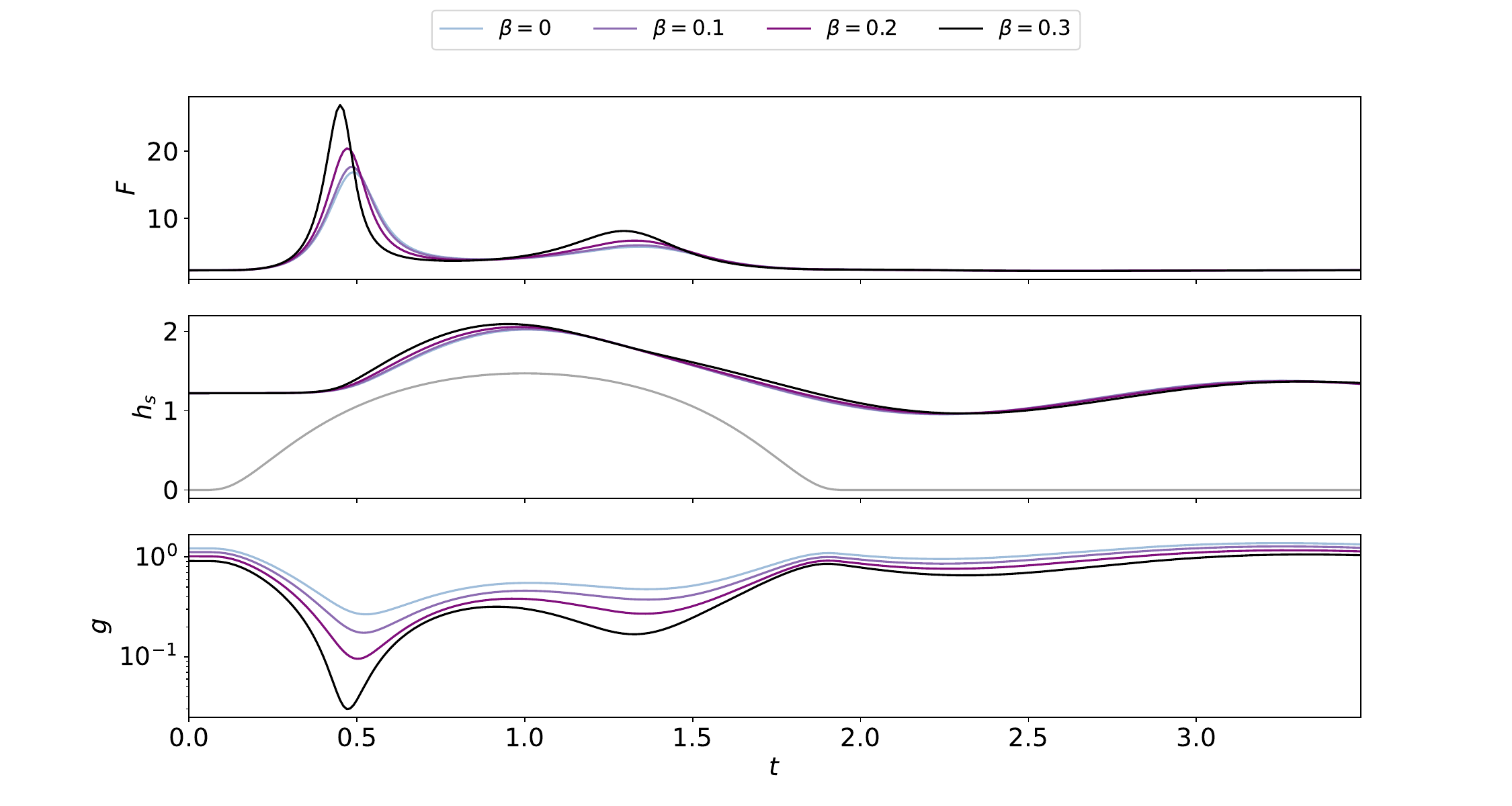}
	\caption{ Example 1: pressure force on stator ($F$, top), stator height ($h_s$,
	middle) and minimum face clearance ($g$, bottom) for a wide bearing ($a = 0.2$),
	showing the behaviour of the bearing when it is disturbed from equilibrium with
	magnitude $\epsilon = 1.2$. Solid grey line in middle plot indicates the axial
	height of the centre of the rotor, $h_R(t)$, and darker colour values represent
	larger misalignment angles. The minimum face clearance is plotted on a logarithmic
	scale. Note that the magnitude of the disturbance is larger than the equilibrium
	clearance between rotor and stator. }
 \label{fig:disturbance_equilibrium}
\end{figure}

A single smooth perturbation to the equilibrium state is modelled by setting
the time-dependent rotor height to be a smooth bump function:

\begin{equation}
	h_R(t)
	=
	\epsilon \exp \left(- \frac{1}{1 - \tfrac 1 4 (t - 2)^2}\right),
\end{equation}
where $\epsilon$ controls the magnitude of the perturbation from equilibrium. The
simulated stator response to this smooth forcing is shown in Figure
\ref{fig:disturbance_equilibrium}. As expected, after the initial disturbance, the
system relaxes to equilibrium. For $t<1$ as the height of the rotor is increasing, peaks
in the pressure force occur earlier and with larger magnitude as the angle of
misalignment increases (top plot). This in turn leads to an earlier and larger
deflection of the stator from its equilibrium position (middle plot). The maximum stator
height does not vary significantly with increasing angle, however there is a significant
decrease in minimum face clearance $g$ (note that $g$ is plotted on a logarithmic
scale). Finally, it is noted that even for significant misalignment angles, a gap is
maintained between the rotor and stator.

\subsection{Example 2: periodic sinusoidal forcing}\label{sec:sine_disturb}

\begin{figure}
	\includegraphics[width=\linewidth]{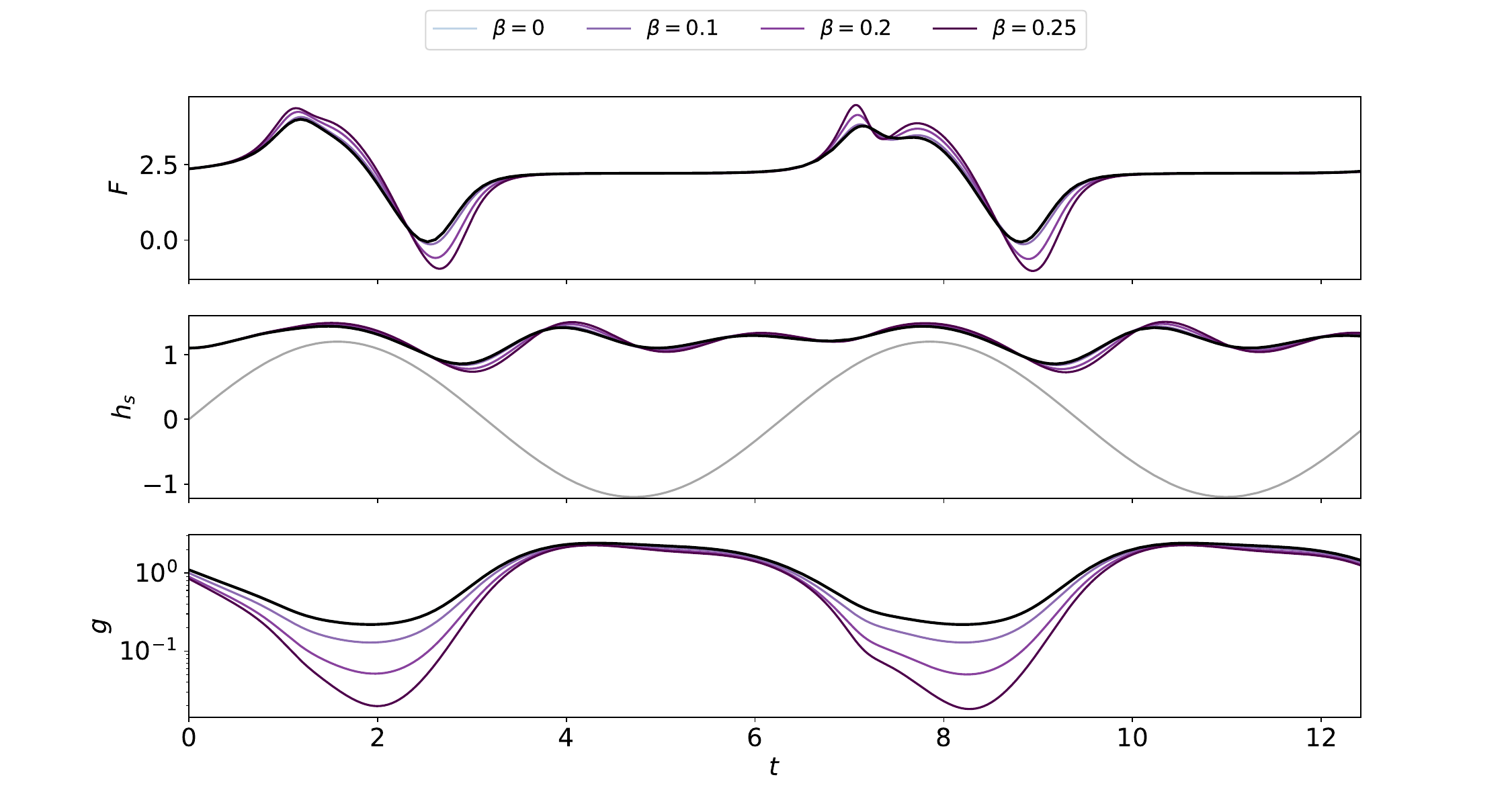}
	\caption{ Pressure force on stator ($F$, top), stator height ($h_s$, middle) and
	minimum face clearance ($g$, bottom) for a wide bearing ($a=0.2$), $\epsilon = 1.2$
	and for various misalignment angles $\beta$.  Solid grey line in middle plot
	indicates the axial height of the centre of the rotor, $h_R(t)$, and darker colour
	values represent larger misalignment angles. Two full periods of the rotor motion is
	shown. For larger misalignment angles, the solution takes between 4 and 5 periods to
	reach a periodic cycle.}
 \label{fig:different_angles_eps_12_mono}
\end{figure}

As a more relevant and challenging test case, the disturbance is now modelled by setting
the time-dependent rotor height as a sinusoidal displacement to model unwanted
vibration:

\begin{equation}
	h_R(t)
	=
	\epsilon \sin(t).
\end{equation}
Thus, the rotor is assumed to move axially.
Stator motion is simulated for a wide bearing ($a = 0.2$), with external
pressurisation, $p_O = 2$, $p_I = 1$, with the ambient pressure $p_a$ taken to
be 1. For this numerical experiment, the amplitude is set as $\epsilon = 1.2$. The
simulation is initialised with a stationary stator so that the test case may be viewed 
as a disturbance to equilibrium. Results of the simulation are shown in
Figure \ref{fig:different_angles_eps_12_mono}.

With increasing misalignment angle, as expected, the minimum face clearance $g$
decreases. The maximal pressure force develops into a larger local peak. The
stator height follows a similar trajectory for all angles, though peaks and
troughs are accentuated with increasing angle. We note that by setting $\beta =
0$, the axisymmetric case is recovered, for which pressure responses can be derived
analytically \cite{bailey2014dynamics}. This serves as a useful benchmark for
verification purposes.

\begin{figure}
	\centering
	\includegraphics[width=0.6\linewidth]{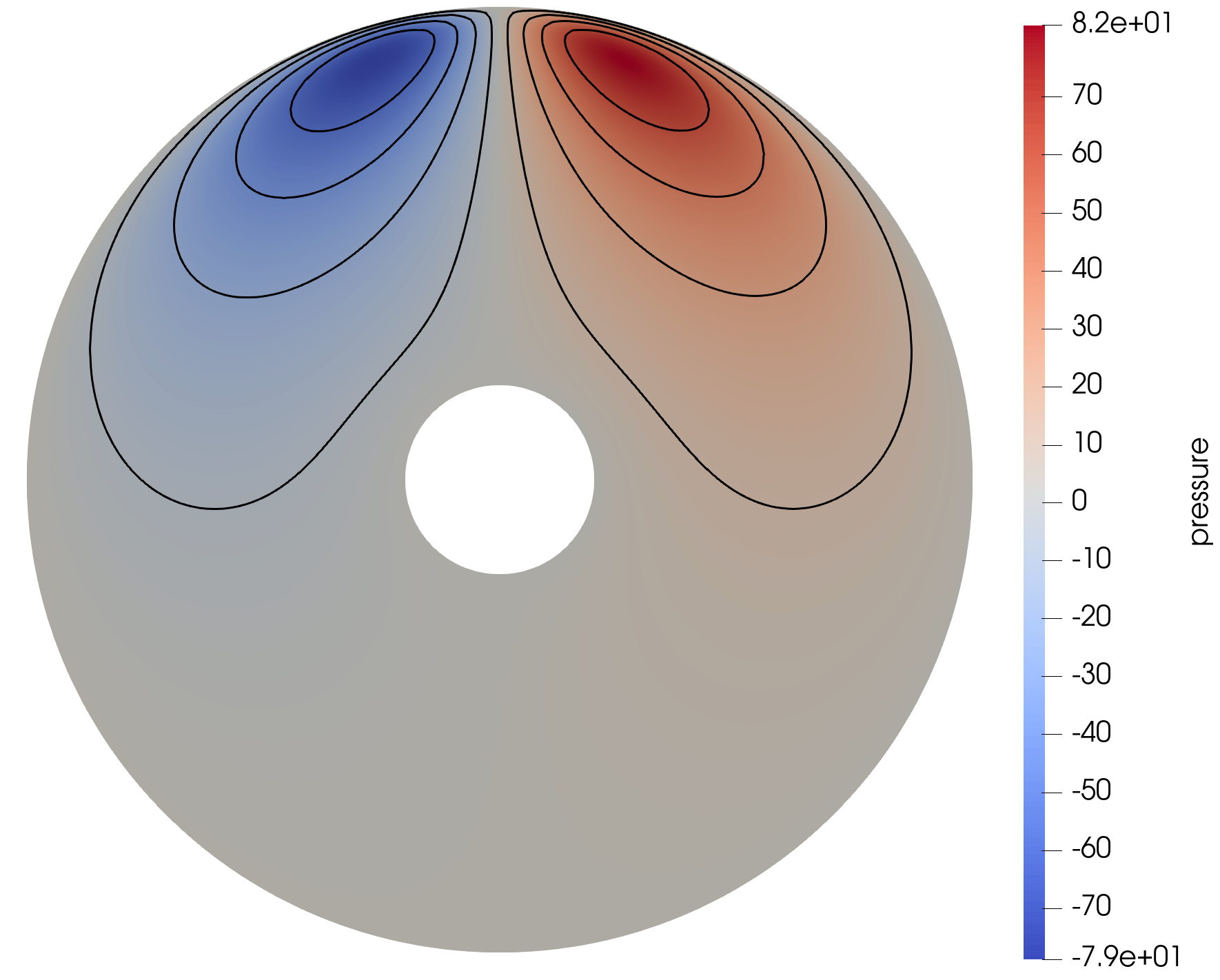}
	\caption{ Contours of the pressure field $p(t, \mathbf x)$ with $t$ chosen such that
	$g(t)$ is at its minimum. Parameter values are $\epsilon = 1.2$, $\beta = 0.25$. The
	rotor motion is anticlockwise. Note the lack of rotational symmetry and steep
	gradients near the outer boundary of the bearing.}
 \label{fig:pressure_misaligned}
\end{figure}

A visualisation of the pressure field at the point at which the minimum face
clearance $g$ is minimal is shown in Figure \ref{fig:pressure_misaligned}. 
There is now significant variation in the azimuthal direction due to
the bearing misalignment, with regions of high and low pressure occurring either
side of the point at which rotor and stator are closest. Indeed, as the minimum
face clearance approaches zero, the pressure field appears to approach a
degenerate solution with singular gradient around the point $(0, 1)$.

\begin{figure}
    \centering
	\includegraphics[width=0.3\linewidth]{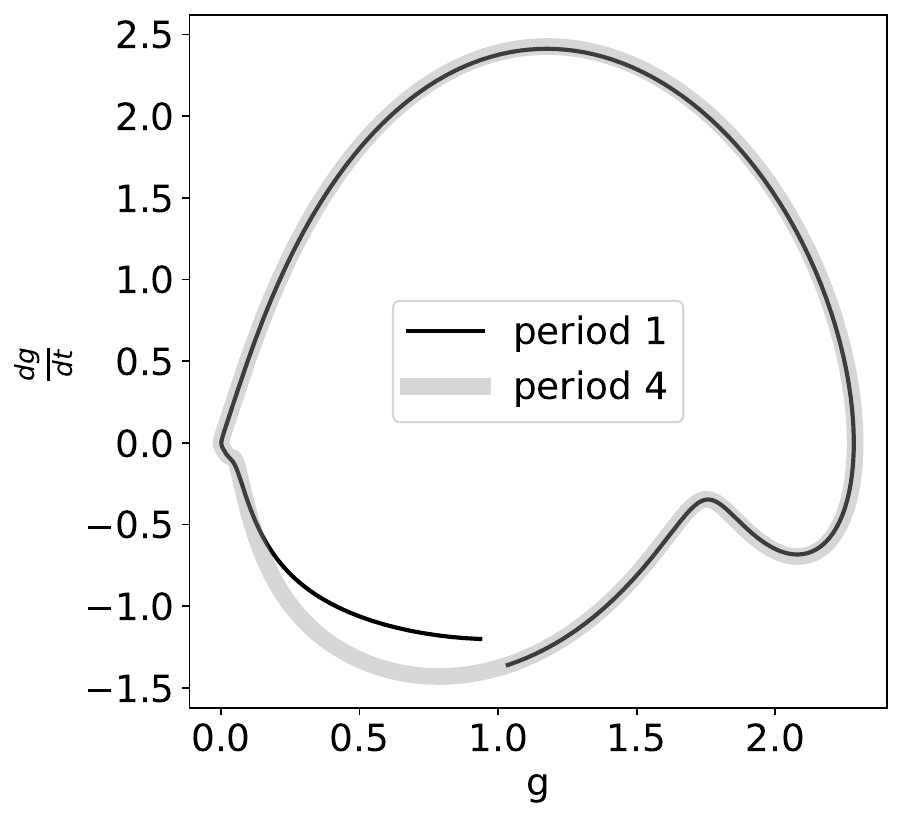}
    \includegraphics[width=0.3\linewidth]{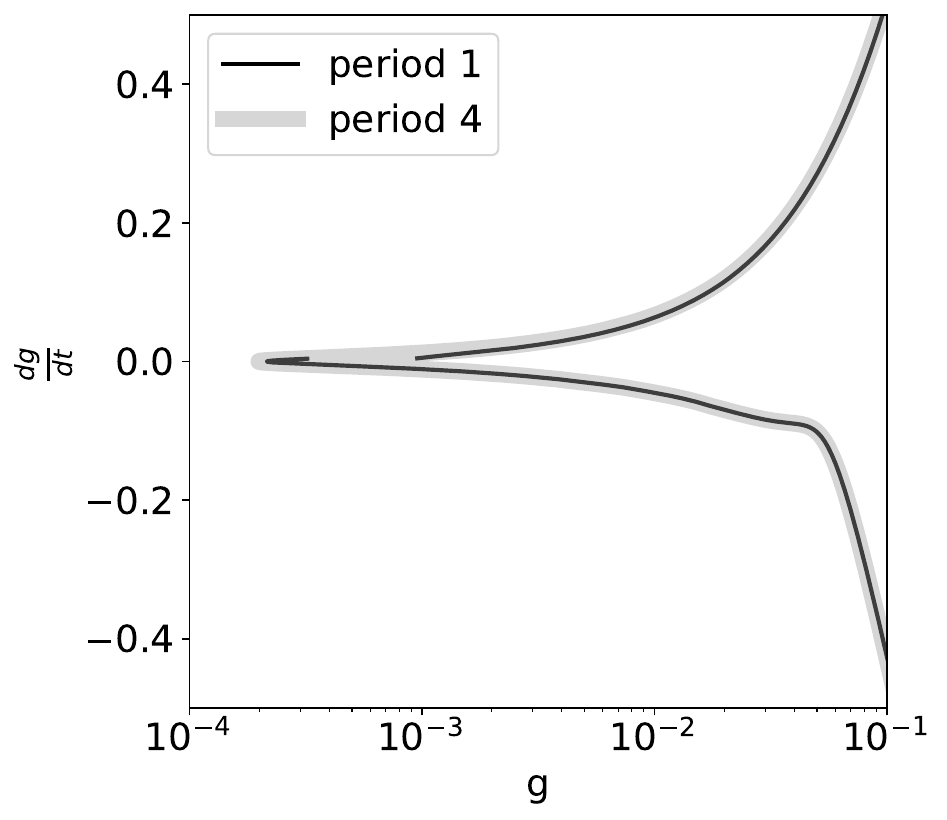}
    \includegraphics[width=0.3\linewidth]{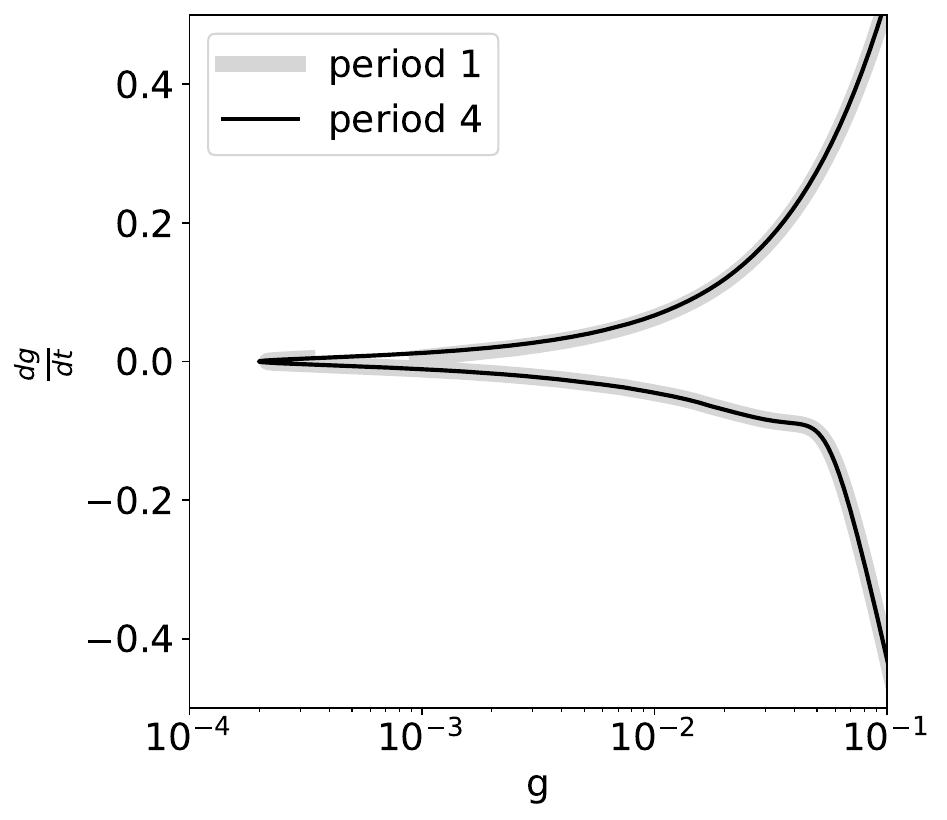}
	\caption{ Illustration of the periodicity of the minimum face clearance, plot of 
    $g$ against $\diff g/\diff t$ for parameter choices $\epsilon = 1.2$ and $\beta = 0.
    2775$, resulting in a very small minimum face clearance. Left: two full periods 
    plotted together, note that the first period starts from equilibrium and the curve 
    is discontinuous. Middle and right: closer view of the portion of the cycle 
    corresponding to the minima of $g$ (i.e. the closest approach of rotor and stator, 
    note also the logarithmic scale) showing that the solution is initially aperiodic 
    (middle), but by the fourth period of the forcing function, periodic dynamics have 
    been reached (right).}
 \label{fig:periodic_or_not}
\end{figure}

It is expected that the minimum face clearance will asymptotically approach a periodic
function in response to the periodic forcing, assuming that rotor-stator contact does
not occur. To investigate this behaviour, $g$ is plotted against $\diff g /\diff t$ over
a period of the forcing function, with a closed curve indicating periodicity. Figure
\ref{fig:periodic_or_not} illustrates that the initial transient phase rapidly
approaches a periodic cycle.

\subsection{Investigation of rotor-stator contact}\label{sec:safety}

The numerical experiments of the preceding sections suggest that for a given
forcing amplitude there exists a critical angle beyond which the pressure force
is not sufficient to maintain separation between rotor and stator. The aim of
this section is to determine the relationship between forcing amplitude and the
associated critical angle. Care must be taken to ensure that results do not
depend upon the time step or the spatial mesh size. When the gap becomes small,
the pressure force becomes large and a small time step will be required. In
addition, the Reynolds equation becomes increasingly convection dominated as the
gap approaches zero which typically results in solutions that are not well
resolved by uniform meshes. We have therefore implemented time step adaptivity
to provide finer resolution when the minimum face clearance becomes small as
well as spatial mesh adaptivity as described in \S\ref{sec:adapt}.

Rotor-stator contact is considered to have occurred when the minimum face clearance
becomes smaller than a small positive tolerance $\texttt{FACE\_TOL}$. This tolerance is
taken to be strictly positive rather than zero due to the fact that the rotor surface
will have small fluctuations and the continuum fluid model begins to break down as the
size of the gap approaches the mean free path of the lubricating fluid. In addition, the
numerical model is not well-defined if the gap becomes zero at any point.

Numerical simulations suggest that the ODE problem \eqref{eqn:stator1order2}
becomes extremely stiff for small face clearances. It is in general not possible to pick
discretisation parameters that can obtain a stable numerical solution for
arbitrarily small gap since the spatial mesh must resolve the layers that occur,
and smaller time steps must be used. However, in light of the above
observations, we are required only to obtain a stable simulation for minimum
face clearances that exceed \texttt{FACE\_TOL}. 

\begin{figure}
	  \includegraphics[width=\linewidth]{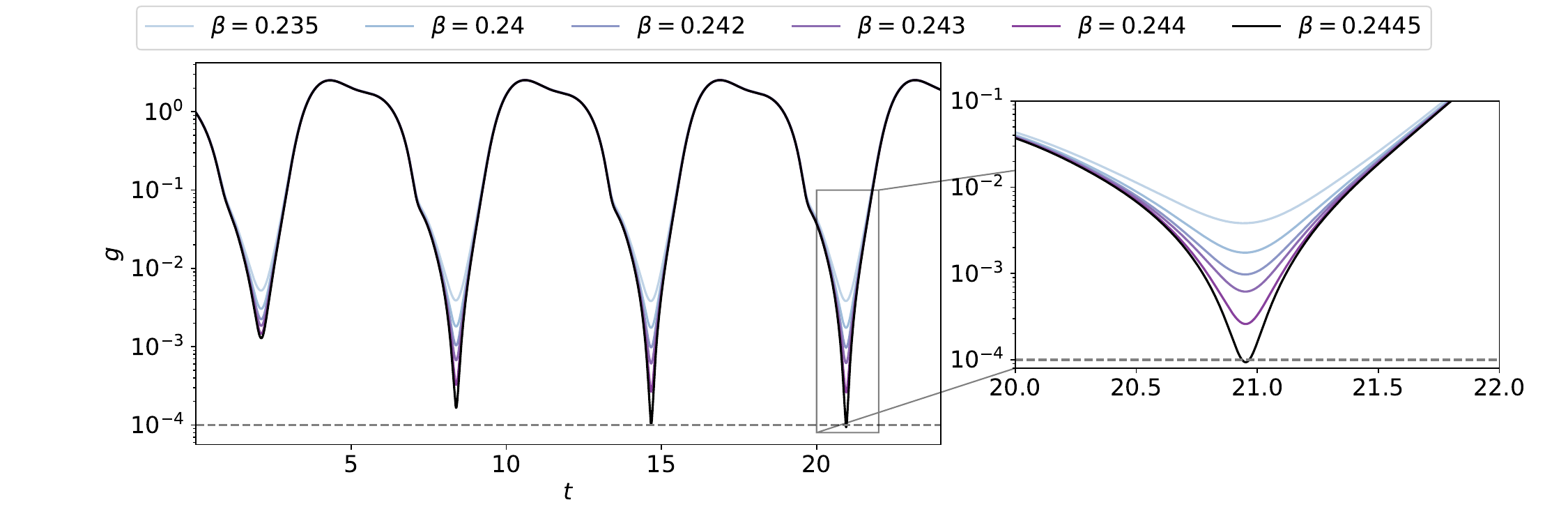}
	  \caption{ Dynamics of the minimum face clearance curve $g$ over a time interval 
      of four periods of the vibrational forcing with amplitude $\epsilon = 1.3$. 
      Dynamics are shown for a selection of misalignment angles $\beta$.  The minimum 
      face clearance is plotted on a logarithmic scale and darker colour values 
      represent larger misalignment angles. The final curve plotted violated the 
      tolerance. To three significant figures the value of the critical angle in this 
      case is 0.244.}
   \label{fig:selection_of_angles}
\end{figure}

\begin{figure}
	\includegraphics[width=0.45\linewidth]{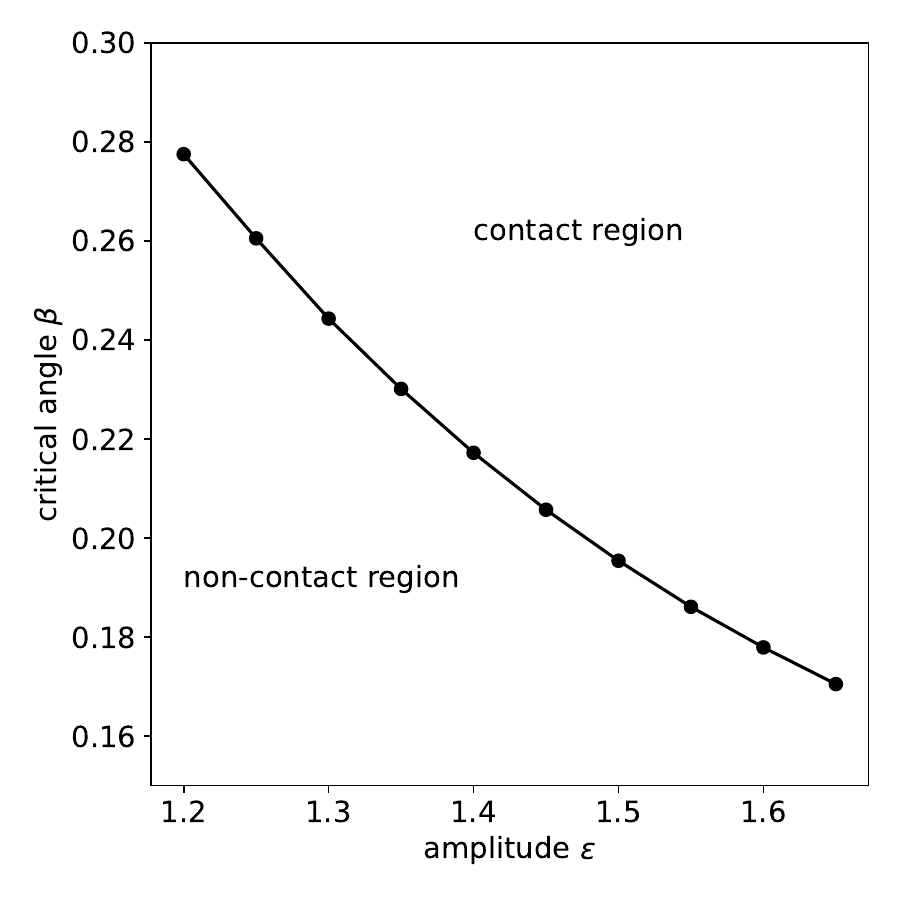}
	\caption{ Plot of the critical misalignment angle against amplitude of disturbance 
    to rotor. Every point $(\varepsilon, \beta)$ lying below the line represents a 
    scenario where the bearing will operate without rotor-stator contact.}
 \label{fig:contact_region}
\end{figure}

For this study, the tolerance is set as \texttt{FACE\_TOL} = $10^{-4}$ and numerical
experiments are conducted to determine the critical misalignment angle for various
forcing amplitudes. Results of the numerical model for $\epsilon = 1.3$ and a variety of
large misalignment angles are plotted in Figure \ref{fig:selection_of_angles}. It is
observed that rather small changes in misalignment angle can result in significant
changes in the minimum clearance. A systematic search was conducted to find, for a given
set of forcing amplitudes $\epsilon_i$, the corresponding largest misalignment angle
$\beta_i$ for which the minimum face clearance stayed above \texttt{FACE\_TOL}. The
pairs $(\epsilon_i, \beta_i)$ can be used to infer a region of safe operation, where the
pressure force is able to withstand disturbances of a given magnitude. The contact and
non-contact regions of the $\epsilon-\beta$ phase space are illustrated in Figure
\ref{fig:contact_region}. 

\begin{figure}
	  \includegraphics[width=\linewidth]
      {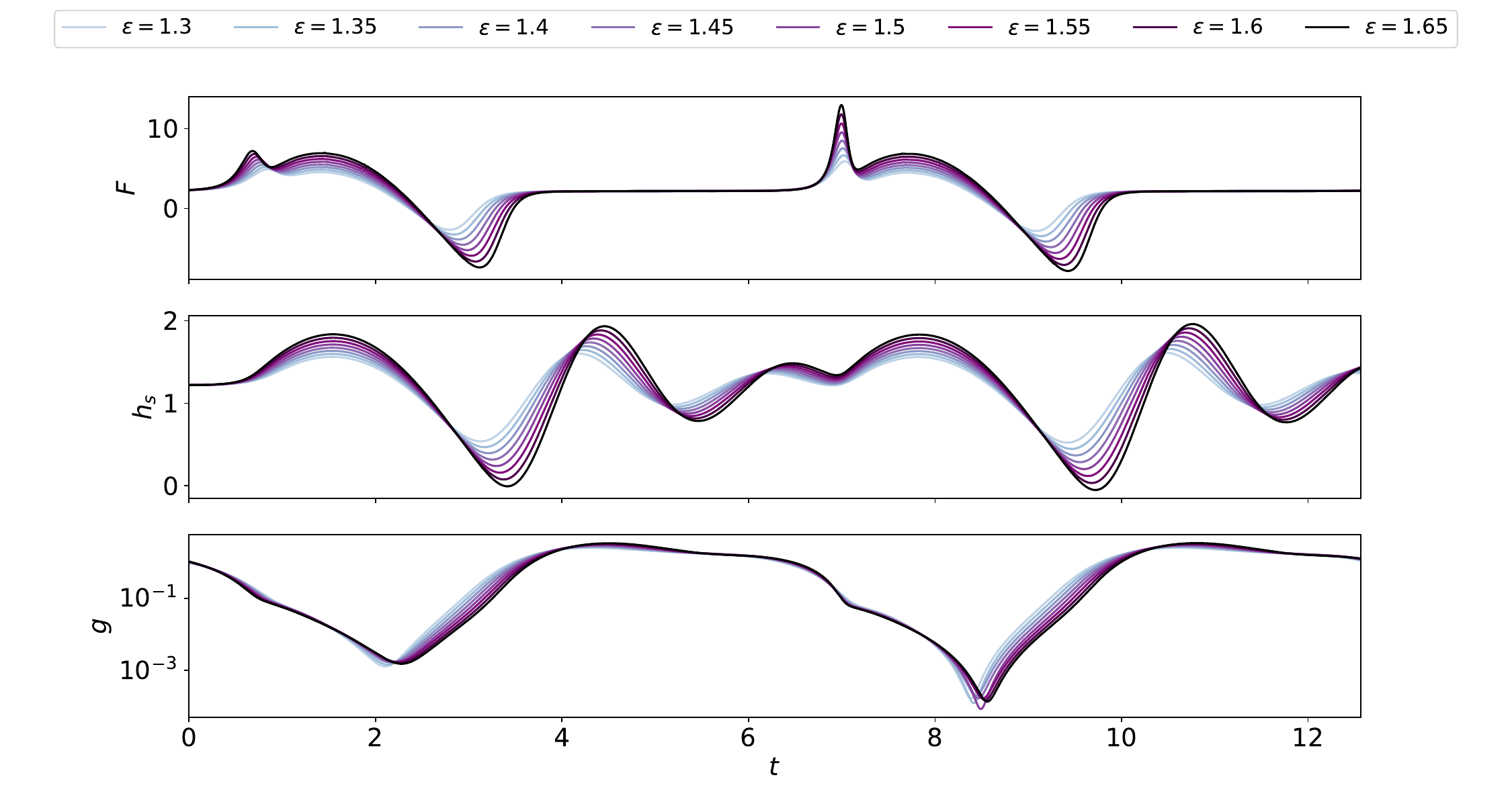} \caption{Pressure force on
	  stator ($F$, top), stator height ($h_s$, middle) and minimum face clearance ($g$,
	  bottom) for a wide bearing ($a=0.2$). The minimum face clearance is plotted on a
	  logarithmic scale and darker colour values represent larger forcing amplitudes.
	  Each line represents the critical misalignment angle for a given forcing
	  amplitude, that is, the smallest misalignment angle such that the minimum face
	  clearance becomes smaller than \texttt{FACE\_TOL}. }
   \label{fig:scritical_angles_vs_amp}
\end{figure}

The dynamics of the system at critical level of misalignment, that is, the largest
$\beta$ for which $g$ remains larger than \texttt{FACE\_TOL}, for a selection of forcing
amplitudes are illustrated in Figure \ref{fig:scritical_angles_vs_amp}. It is noted that
the minimum gap is not exactly \texttt{FACE\_TOL} in all cases. This is due to the fact
that the critical angle was determined to 4 decimal places, but the minimum gap is very
sensitive to changes in misalignment angle 

A natural extension of the present work would be to include more appropriate
boundary conditions for the case where the face separation is very small. A
Navier slip boundary condition is considered for example in
\cite{bailey2015dynamics} for the case of a parallel bearing and
\cite{bailey2016dynamics} for a coned bearing.

\section{Conclusions}

A coupled mathematical model of a misaligned non-contacting mechanical face seal was
presented, where the fluid flow was based on a lubrication approximation of the
Navier-Stokes equations and the stator as a spring-mass-damper system. The governing
equations were solved using a developed numerical model, utilising adaptive meshing
techniques to provide fast and accurate numerical approximations to the governing
equations. Numerical experiments were conducted and a safe region of operation of the
seal was determined, in terms of misalignment of rotor and stator and the magnitude of
external disturbances to equilibrium.

In principle the presented approach may be extended to arbitrary rotor and stator
geometries without difficulty, allowing different designs or manufacturing defects to be
tested for safety in a variety of undesirable operating conditions, and safe operating
parameter sets determined.

\section*{Acknowledgements}

This work originated from a seed grant funded by the Institute for Mathematical
Innovation at the University of Bath. All authors were supported by the EPSRC grant
EP/X030067/1. BA and TP received support from the Leverhulme Trust grant RPG-2021-238.
TP was also supported by the EPSRC programme grant EP/W026899/1.

\appendix

\section{Derivation of Misaligned Boundary Conditions}\label{sec:bcs_stuff}

In this section details are provided on the rotated boundary conditions resulting from
misalignment of the bearing given in \S\ref{subsec:airflowmodel}. We assume that only
the rotational component is altered, and that the axial movement of the rotor remains
parallel to the $\hat{z}$-axis. In the aligned case the boundary condition at the rotor
is given by 
\begin{equation}
    \mathbf{\hat u}(\hat x, \hat y, \hat z) 
    =
    (- \hat \Omega \hat y, \hat \Omega \hat x, 0)^T,
\end{equation}
in Cartesian coordinates. Let
\begin{equation}
  R =
    \begin{pmatrix}
        1 & 0 & 0 \\
        0 & \cos \hat \beta & - \sin \hat \beta \\
        0 & \sin \hat \beta & \cos \hat \beta
    \end{pmatrix}
\end{equation}
be the matrix representing a rotation by angle $\hat \beta$ around the $\hat x$-axis in
the standard basis of $\mathbb{R}^3$. Then at a point $\mathbf{\hat x}$ on the
misaligned rotor, the velocity boundary condition is given by
\begin{equation}
    \mathbf{\hat u}_{\hat \beta}(\mathbf{ \hat x})
    =
    R \mathbf{\hat u}(R^{-1} \mathbf{\hat x}).
\end{equation}
After transforming to cylindrical coordinates, one obtains
\begin{equation}
      \hat{u}_{\hat \beta} = - \hat{\Omega} \hat{z} \sin \hat{\beta}\cos \hat{\theta}, 
      \quad 
      \hat{v}_{\hat \beta} 
      = 
      \hat{\Omega} \hat{r} \cos \hat{\beta} 
      + 
      \hat{\Omega} \hat{z} \sin \hat{\beta} \sin \hat{\theta}, \quad 
      \hat{w}_{\hat \beta} = \hat{\Omega} \hat{r} \sin \hat{\beta} \cos \hat{\theta} 
      + 
      \frac{\partial \hat{h}_r}{\partial \hat{t}},
      \text{  at  } \hat z = \hat{h}_r,
\end{equation}
from which \autoref{eq:bc_rotor_dimensional} follows after applying small angle
approximations.

\subsection{Dimensional Analysis}\label{sec:appendix_dimensional}

In this section we derive the non-dimensional velocity boundary conditions given in
\autoref{eq:noslipnondimvelbc}. Using the same scalings as in \S\ref{subsubsec:dim}, the
above become
\begin{equation}
    u_{\beta}
    =
    - \frac{\hat{\Omega} \hat{h}_0 \hat{\delta}_0}{\hat{U}} \beta z \cos \theta
    =
    - 
    \frac{\hat{r}_0 }{\hat{r}_0} 
    \frac{\hat{\Omega} \hat{h}_0 \hat{\delta}_0}{\hat{U}} \beta z \cos \theta
    =
    -\frac{Re^* \hat{h}_0 \hat{\delta}_0}{\hat{\delta}_0 \hat{r}_0}\beta z \cos \theta
    =
    O(Re^* \hat{\delta}_0)
\end{equation}
\begin{equation}
    v_{\beta}
    =
    \frac{\hat \Omega \hat{r}_0}{\hat V}r
    +
    \frac{\hat \Omega \hat{h}_0 \hat{\delta}_0}{\hat V} \beta z \sin \theta
    =
    r 
    +
    O(\hat{\delta}_0^2)
\end{equation}
\begin{equation}\label{eq:z_bc_scaling}
    w_{\beta}
    =
    \frac{\hat{T}\hat{\Omega}\hat{r}_0\hat{\delta}_0}{\hat{h}_0} \beta \cos \theta
    +
    \frac{\partial h_r}{\partial t}.
\end{equation}
Examining the size of the first term in \autoref{eq:z_bc_scaling}, we calculate, using
the definition of the squeeze number in \autoref{eqn:deltasigmaFr},
\begin{equation}
    \frac{\hat{T}\hat{\Omega}\hat{r}_0\hat{\delta}_0}{\hat{h}_0}
    =
    \hat{\Omega} \hat{T}
    =
    \frac{\hat{\Omega} \hat{r}_0}{\hat{U} \sigma}
    =
    \frac{Re^*}{\sigma},
\end{equation}
indicating that this term must be considered part of the leading order velocity. Finally
then the axial component of the boundary condition is 
\begin{equation}
    w_{\beta}
    =
    \frac{\partial h_r}{\partial t}
    +
    \frac{Re^*}{\sigma}r \beta \cos \theta.
\end{equation}

\section{Solution of the Reynolds Equation by the Finite Element Method}
\label{sec:pde_stuff}

In this section, further details are provided on the solution of the Reynolds equation
by the finite element method, beginning with the formulation of the Reynolds equation in
an appropriate functional setting and discussion of existence \& uniqueness of solutions
in \S\ref{sec:exist}, followed by the definition of the finite element approximation to
the pressure field and description of the calculation of the forcing $\mathcal{F}$ in
\S\ref{sec:fem_appendix}. The procedure for dynamic mesh adaptivity is described in
\S\ref{sec:meshing_appendix}.

\subsection{Existence \& Uniqueness of Solutions of the Reynolds Equation}
\label{sec:exist}

The modified Reynolds equation \eqref{eq:final_Reynolds}, stated again here for
convenience, 

\begin{equation}\label{eq:cartesian_Reynolds7}
   \tilde{\sigma}\frac{\partial}{\partial t}h 
   - 
   \nabla \cdot \left(h^3 \nabla p\right)
   =
   6 Re^*r \beta \cos \theta,
\end{equation}
is a time-dependent diffusion equation. The nature of the nonlinear diffusion
coefficient means that the problem approaches degeneracy in the limit $h \to 0$, which
can lead to numerical instability when the spatial resolution is not sufficient. We opt
to use a standard Galerkin finite element method with adaptive meshing to ensure
sufficient resolution of the problem features. This is motivated by the fact that steep
pressure layers will occur where $h$ is smallest, and will be highly localised in space
and time, meaning that computation time can be reduced by using coarse meshes when and
where possible. The reader is referred to
\cite{ainsworth1997posteriori,bangerth2003adaptive,nochetto2009theory} for general
references on the rationale for and techniques associated with adaptive finite element
methods.

The finite element method for this problem is based upon the weak formulation of a time
semidiscrete version of \eqref{eq:cartesian_Reynolds7}. To provide the proper functional
setting, let $\operatorname{L}^{2}(\Omega)$ denote the Lebesgue space of square Lebesgue
integrable functions over $\Omega$, $\operatorname{H}^{1}(\Omega)$ denote the Hilbert
space of functions in $\operatorname{L}^{2}(\Omega)$ whose first (weak) derivatives are
in $\operatorname{L}^{2}(\Omega)$, and $\operatorname{H}^{1}_0(\Omega)$ be the subspace
of $\operatorname{H}^{1}(\Omega)$ consisting of functions which vanish on the boundary
$\partial \Omega$. Weak solutions of \autoref{eq:cartesian_Reynolds7} are sought in the
function space

\begin{equation}
    V := \{q \in \operatorname{H}^1(\Omega): 
    q = p_I \,\,\text{on} \,\, \Gamma_I, \,\, q=p_O \,\, \text{on}\,\, \Gamma_O \}.
\end{equation}
The weak formulation is then to find $p_*\in V$ such that 
\begin{equation}\label{eq:weak_form}
	\int_{\Omega}
	h_*^3 \nabla p_* \cdot \nabla q 
	=
	\int_{\Omega}
    fq \diff x
	\quad \forall q \in H^1_0(\Omega),
\end{equation}
where 
\begin{equation}
    f = 6 Re^*r \beta \cos \theta - \tilde{\sigma}\frac{\partial h_*}{\partial t}.
\end{equation}
We make the assumption that there is a constant $h_{\text{min}}>0$ such that $h_* \geq
h_{\text{min}}$ for all choices of $h_{s,\,*}$ and $h_{r,\,*}$ made henceforth. This
amounts to assuming that the fluid film has strictly positive thickness bounded away
from zero for all time, which is to say that the rotor and stator do not come into
contact. Under this assumption, standard results on elliptic PDEs (see for example
Theorem 3.8 of \cite{ern2004theory}) may be applied to conclude that problem
\ref{eq:weak_form} has a unique solution.

\begin{remark}[Alternative bearing geometries]
    A popular design of bearings includes troughs in the bearing surface in a spiral
    pattern. A parameterisation of such a bearing, $h(r, \theta)$ would be a step
    function and therefore not differentiable. However in the theory of elliptic PDEs,
    the diffusion coefficient is required only to be essentially bounded (that is,
    contained in $\operatorname{L}^{\infty}(\Omega)$), meaning that such designs can be
    simulated in this framework.
\end{remark}

\subsection{The finite element method}\label{sec:fem_appendix} 

To apply the finite element method, the annular domain $\Omega$ is approximated by a
domain $\Omega'$ for which the inner and outer boundaries are polygonal approximations
to the circular geometry. The approximate domain $\Omega'$ is subdivided into a
conforming triangulation of quadrilateral elements, with this subdivision denoted by
$\mathcal{T}$, that is, $\mathcal{T}$ is a finite family of sets such that  
\begin{enumerate}
    \item 
    Every element $K \in \mathcal T$ is an open quadrilateral,
    \item 
    For any elements $K,J \in \mathcal{T}$ we have that $\overline K \cap \overline J$
    is either the empty set, a single point (i.e. an element vertex), a common edge of
    $K$ and $J$ or the whole of $\overline K$ and $\overline J$,
    \item 
    $\bigcup_{K\in\mathcal T}\overline K=\overline {\Omega}'$.
\end{enumerate}
We define the approximation spaces 
\begin{equation}\label{eq:fem_spaces}
	\begin{split}
	\mathbb{V} 
	&:=
	\{\chi \in \operatorname{H}^1(\Omega') : \chi \text{ is bilinear on each } 
	T \in \mathcal{T},\, \text{satisfies B.C.s in \eqref{eqn:bndry}}\},\\
	\mathbb{V}_{0} 
	&:=
	\{\chi \in \operatorname{H}^1(\Omega') : \chi \text{ is bilinear on each } 
	T \in \mathcal{T}, \, \chi = 0 \,\,\text{on}\,\, \partial \Omega'\}.
	\end{split}
\end{equation}
The finite element method is then to find the approximate pressure $\Phi \in \mathbb{V}$
such that 
\begin{equation}\label{eq:fem}
	\int_{\Omega'}
	h_*^3\nabla \Phi \cdot \nabla \chi
	\diff x
	=
	\int_{\Omega'}
    \tilde{f} \chi \diff x
	\quad \forall \chi \in \mathbb{V}_0,
\end{equation}
where $\tilde{f}$ is an approximation of $f$, which is included to emphasise that in
practice the integrals in \autoref{eq:fem} are approximately assembled using quadrature.
Existence and uniqueness of the discrete solution follows in a similar manner to the
above since a conforming finite element method is used. 

Choosing bases for the approximation spaces in \autoref{eq:fem_spaces} allows
\autoref{eq:fem} to be assembled to obtain a linear algebraic system, which can be
solved to obtain $\Phi$. Following this calculation, $\Phi$ can be used to obtain the
approximate leading order pressure force on the stator \eqref{eq:leading_order_force},
which is denoted by $\mathcal{F}$, that is,

\begin{equation}\label{eq:appendix_approx_F}
	\mathcal{F}\left(h_{s,\,*}, h_{r,\,*}, \frac{\diff h_*}{\diff t}\right)
	:= 
	\int_{\Omega'}(\Phi - p_a) \diff x \diff y.
\end{equation}

\subsection{Mesh Adaptivity}\label{sec:meshing_appendix}

In this work, dynamic mesh adaptivity is used to minimise computational cost while
providing fine local resolution is regions of interest. Here, further details of the
adaptive algorithm are provided.

\subsubsection{Local Error Indicators}
The Kelly error indicator \cite{kelly1983posteriori} is chosen as a criterion by which
to decide which cells to refine or coarsen. Piecewise linear and continuous finite
element functions have discontinuous gradient across element edges. The Kelly error
indicator is based on the size of these jumps. These ideas are made precise with the
following definitions.

\begin{definition}[Jump operator]
The jump of a vector quantity $\bm{v}$ over a common edge of elements $K_1$ and $K_2$
with outward facing normals $\bm{n_{K_1}}$ and $\bm{n_{K_2}}$ is defined as
\begin{equation}
    \left\llbracket \bm{v} \right\rrbracket
    :=
    \bm{v}_{K_1} \cdot \bm{n}_{K_1} + \bm{v}_{K_2} \cdot \bm{n}_{K_2}.
\end{equation}
\end{definition}

\begin{definition}[Error indicator]
For an element $K \in \mathcal{T}$, the local error indicator $\eta_K$ is defined by
\begin{equation}\label{eq:kelly}
    \eta_K(\Phi)^2
    :=
    \operatorname{diam}(K)\int_{\partial K}
    \left\llbracket \nabla \Phi \right\rrbracket^2.
\end{equation}
\end{definition}
The Kelly error indicator has been shown to be the dominant part of the finite element
approximation error for Poisson's equation through rigorous error analysis
\cite{ainsworth1997posteriori}, but is often applied more widely as a heuristic but
effective criterion for mesh refinement. Its computation is relatively cheap as it
requires only numerical integration of the discrete solution.

\subsubsection{Implementation Details}\label{sec:implement_adaptivity}

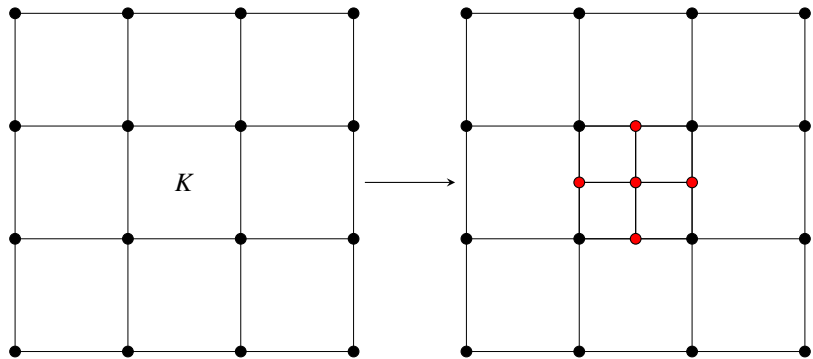
\begin{figure}[h!]
    \begin{center}
        \begin{tikzpicture}[scale=0.75]
        \draw[step=2cm,black,very thin] (0,0) grid (6,6);
        \draw[black,fill=black] (0,0) circle (.6ex);
        \draw[black,fill=black] (0,2) circle (.6ex);
        \draw[black,fill=black] (0,4) circle (.6ex);
        \draw[black,fill=black] (0,6) circle (.6ex);
        \draw[black,fill=black] (2,0) circle (.6ex);
        \draw[black,fill=black] (2,2) circle (.6ex);
        \draw[black,fill=black] (2,4) circle (.6ex);
        \draw[black,fill=black] (2,6) circle (.6ex);
        \draw[black,fill=black] (4,0) circle (.6ex);
        \draw[black,fill=black] (4,2) circle (.6ex);
        \draw[black,fill=black] (4,4) circle (.6ex);
        \draw[black,fill=black] (4,6) circle (.6ex);
        \draw[black,fill=black] (6,0) circle (.6ex);
        \draw[black,fill=black] (6,2) circle (.6ex);
        \draw[black,fill=black] (6,4) circle (.6ex);
        \draw[black,fill=black] (6,6) circle (.6ex);
        \draw (3,3) node {$K$};
        \draw [-stealth](6.2,3) -- (7.8,3);
        \draw[step=2cm,black,very thin] (8,0) grid (14,6);
        \draw[color=black] (10,2) rectangle (11,3);
        \draw[color=black] (10,3) rectangle (11,4);
        \draw[color=black] (11,3) rectangle (12,4);
        \draw[color=black] (11,2) rectangle (12,3);
        \draw[black,fill=red] (11,2) circle (.6ex);
        \draw[black,fill=red] (11,4) circle (.6ex);
        \draw[black,fill=red] (10,3) circle (.6ex);
        \draw[black,fill=red] (12,3) circle (.6ex);
        \draw[black,fill=red] (11,3) circle (.6ex);
        \draw[black,fill=black] (8,0) circle (.6ex);
        \draw[black,fill=black] (8,2) circle (.6ex);
        \draw[black,fill=black] (8,4) circle (.6ex);
        \draw[black,fill=black] (8,6) circle (.6ex);
        \draw[black,fill=black] (10,0) circle (.6ex);
        \draw[black,fill=black] (10,2) circle (.6ex);
        \draw[black,fill=black] (10,4) circle (.6ex);
        \draw[black,fill=black] (10,6) circle (.6ex);
        \draw[black,fill=black] (12,0) circle (.6ex);
        \draw[black,fill=black] (12,2) circle (.6ex);
        \draw[black,fill=black] (12,4) circle (.6ex);
        \draw[black,fill=black] (12,6) circle (.6ex);
        \draw[black,fill=black] (14,0) circle (.6ex);
        \draw[black,fill=black] (14,2) circle (.6ex);
        \draw[black,fill=black] (14,4) circle (.6ex);
        \draw[black,fill=black] (14,6) circle (.6ex);
        \end{tikzpicture}
    \end{center}
    \caption{Adapting a mesh by refining the central element $K$. This results in new 
    element vertices being created, shown in red.}
    \label{fig:refinement}
\end{figure}

At time $t=0$, a coarse mesh is generated with 320 elements and $1280$ degrees of
freedom. This coarse mesh is sufficient to be used when the rotor and stator are
relatively far apart, with adaptive mesh refinement activated when the minimum face
clearance is small. In this work, mesh adaptivity was used when $g(t^n) < 2 \times
10^{-2}$. Every \texttt{refine\_interval} time steps, the local error indicator is
calculated on each element and used to determine whether an element is refined,
coarsened, or left untouched. Pseudocode is provided in \autoref{alg:adaptivity}. The
mesh obtained at the end of the algorithm is the initial mesh for the following time
step. Quadrilateral elements are refined by quadrisection where the midpoint of each
side is joined to the element centroid to make a new edge (see
\autoref{fig:refinement}). Thus, each refined element has 4 children. Elements are only
coarsened if all four \lq child' cells of a larger element that has been previously
refined are marked for coarsening. It should be noted that certain limits are applied to
the refinement and coarsening. Cells cannot be coarsened beyond the initial coarse mesh
level, and a limit on the number of refinement levels is set to control the amount of
computational resources required.

\begin{algorithm}
    \SetKwInOut{Input}{Input}
    \SetKwInOut{Output}{Output}

    Given $h_{s,\,*}$, $h_{r,\,*}$, $\diff h_* \slash \diff t$, $\theta_{\text{max}}$
    and $\theta_{\text{min}}$: \\
    Solve \autoref{eq:fem} to find $\Phi$;\\
    \If{$l \equiv 0$\, ($\mod$ \texttt{refine\_interval})}
    {
    \If{$g(t) < 2.0\times 10^{-2}$}
      {
    \For{$K \in \mathcal T$}
    {
            Compute $\eta_K(\Phi)$;\\
            \If{$\eta_K(\Phi) > \theta_{\text{max}} \max_K \eta_K(\Phi)$}
            {Mark $K$ for refinement;}
            \ElseIf{$\eta_K(\Phi) < \theta_{\text{min}} \min_K \eta_K(\Phi)$}
            {Mark $K$ for coarsening;}
        }
        Apply refinement limits and perform refinement \& coarsening.
      }
      \Else
      {
        Perform coarsening step.
      }
      }
    \caption{Dynamic mesh refinement at timestep $t^l$}
    \label{alg:adaptivity}
\end{algorithm}

\section{Time stepping}\label{sec:time_stepper}

We shall use $\mathcal{F}$ (defined in \autoref{eq:approx_F}, see also
\autoref{eq:appendix_approx_F}) as a numerical approximation of $F$, which is available
for a given time, stator height and rate of change. In practice this will depend on the
prescribed rotor height, and approximations to $x_1$ and $x_2$. Let $0 = t^0 < t^1 <
\dots < t^N = T$ be a partition of the interval $[0, T]$. We denote approximations to
$x_i(t^n)$ by $X_{i}^n$. Initial conditions are given by \eqref{eqn:ntl}. Given
$X_{1}^n$ and $X_{2}^n$, with $\Delta t^n := t^{n+1} - t^n$, numerical approximations
$X_{1}^{n+1}$ and $X_{2}^{n+1}$ are obtained using the RK4 method. For further details
on the RK4 method as well as solvers for ordinary differential equations more broadly,
we refer the reader to \cite{butcher2016numerical} and the references therein. Letting
$\mathbf{X}^n$ denote the vector $(X_{1}^{n},\,X_{2}^{n})$, this method can be expressed
in a compact form. Let the predictor steps $\mathbf{Y}^{n,\,i}$ be

\begin{equation}
    \mathbf{Y}^{n,\,i}
    =
    \mathbf{X}^n + \Delta t^n \sum_{j=1}^4 a_{ij}\mathbf{F}\left(t^{n} 
    + c_j \Delta t^n,\mathbf{Y}^{n,\,i}\right),
\end{equation}
for $i=1,...,4$.

where the components of the vector $\mathbf{F}(t, \mathbf{Y})$ are given by
\begin{equation}
\begin{split}
    \mathbf{F}(t, \mathbf{Y})_1
    &:=
    Y_2,\\
    \mathbf{F}(t, \mathbf{Y})_2
    &:=
    \alpha \mathcal{F}\left(Y_1, h_r(t), Y_2 - \dot{h}_r(t))\right)
    - D_{\alpha} Y_2 - K_z(Y_1 -1).
\end{split}
\end{equation}
Then 
\begin{equation}
    \mathbf{X}^{n+1}
    =
    \mathbf{X}^n 
    +
    \Delta t^n \sum_{i=1}^4 b_i \mathbf{F}
    \left(t^n + c_i \Delta t^n, \mathbf{Y}^{n,\,i}\right),
\end{equation}
where the weights $b_i$ and nodes $c_i$ are given by the Butcher tableau in
\autoref{tab:rk4}.

\begin{table}[h!]
\centering
\renewcommand{\arraystretch}{1.5} 
\setlength{\arraycolsep}{10pt}     

\[
\begin{array}{c|cccc}
0 &  &  &  &  \\
\frac{1}{2} & \frac{1}{2} &  &  &  \\
\frac{1}{2} & 0 & \frac{1}{2} &  &  \\
1 & 0 & 0 & 1 &  \\
\hline
  & \frac{1}{6} & \frac{1}{3} & \frac{1}{3} & \frac{1}{6}
\end{array}
\]

\caption{Butcher tableau for the RK4 method.}
\label{tab:rk4}
\end{table}
 
\bibliographystyle{unsrt}
\bibliography{bibliography_file} 

\end{document}